%% file: hhf_arch.tex
\documentclass[a4paper,12pt]{article}%{pcms2-l}
\usepackage{theorem,amsmath,latexsym,amssymb,newcom,esint}
\begin{document}
\title{\sf{On a new harmonic heat flow \\
           with the reverse H\"older inequalities }}
\author {\sf{Kazuhiro HORIHATA}}
\maketitle
\begin{abstract}
This paper first proposes a new approximate scheme 
 to construct a harmonic heat flow $u$
  between $(0,T)$ $\times$ $\mathbb{B}^d$
   and $\mathbb{S}^D$ $\subset$ $\mathbb{R}^{D+1}$ 
    with positive integers $d$ and $D$:
We agree that harmonic heat flow $u$ means a solution of
\begin{equation}
\frac{\partial u}{\partial t}
 \, - \, \triangle u
  \, - \, | \nabla u|^2 u
   \; = \;  0.
\notag
\end{equation}
It's scheme is crucially given by
\begin{equation}
\frac{\partial u_\lambda}{\partial t}
 \, - \, \triangle u_\lambda
  \, + \, \lambda^{1-\kappa}
   \bigl( | u_\lambda |^2 \, - \, 1 \bigr) u_\lambda \; = \; 0,
\notag
\end{equation}
where the unknown mapping $u_\lambda$ is 
 from $(0,T)$ $\times$ $\mathbb{B}^{d}$
  to $\mathbb{R}^{D+1}$
   with positive numbers $\lambda$ and $T$, 
    positive integers $d$ and $D$
     and $\kappa (t)$ $=$
      $\arctan (t)/\pi$ $(0 \le t)$.
\par
The benefits to introduce a time-dependent parameter
 $\lambda^{1-\kappa}$ is readily to see
\begin{equation}
\int_{Q(T)} \lambda^{1-\kappa}
 ( | u_\lambda |^2 \, - \, 1 )^2 \, dz 
  \; \le \;
   C/\log \lambda
\notag
\end{equation}
for some positive constant $C$ 
 independent of $\lambda$.
\par
Next, making the best of it, we establish that a passing to the limits $\lambda \nearrow \infty$
 ( modulo subsequence of $\lambda$ )
   brings the existence of a harmonic heat flow into spheres with
\renewcommand{\theenumi}{\roman{enumi}}
\begin{enumerate}
\item
a global energy inequality,
\item
a monotonicity for the scaled energy,
\item
a reverse Poincar\'e inequality.
\end{enumerate}
These inequalities (\romannumeral 1), (\romannumeral 2) and (\romannumeral 3) improves
 the estimates on it's singular set
  contrast to the results by Y.~Chen and M.~Struwe \cite{chen-struwe},
i.e. I show that a singular set of the new harmonic heat flow into spheres
 has at most finite $(d-\epsilon_0)$-dimensional Hausdorff measure with respect to the parabolic metric
  whereupon $\epsilon_0$ is a small positive number.
\par
I believe that inequalities (\romannumeral 1), (\romannumeral 2) and (\romannumeral 3)
 allow us to analyze how it behaves around its singularities.
\end{abstract}
%
\input intro_hhf_arch.tex
\input glhf_hhf_arch.tex
\input estimate_hhf_arch.tex
\input hhf_hhf_arch.tex
\input main_hhf_arch.tex

\input ref_hhf_arch.tex
\end{document}

%% file: intro_hhf_arch.tex
%
%#! platex hhf
%
\setcounter{chapternumber}{1}\setcounter{equation}{0}
\renewcommand{\theequation}%
           {\thechapternumber.\arabic{equation}}
\section{\enspace Introduction.}
\label{SEC:Intro}
Let $\mathbb{B}^d$ and $\mathbb{S}^D$ be respectively the unit ball centred
 at the origin in $\mathbb{R}^d$ and
  the unit sphere in $\mathbb{R}^{D+1}$,
   where $d$ and $D$ are positive integers greater than or equal to $2$.
Consider the Sobolev space 
 $H^{1,2} (\mathbb{B}^d ;\mathbb{S}^D)$ $\, := \,$
  $\{u \in H^{1,2} (\mathbb{B}^d ; \mathbb{R}^{D+1})$
    $; |u| = 1 \quad \mathrm{a.e.} \, \, x \in \mathbb{B}^d \}$.
     Giving a mapping $u_0 \, \in \, $
      $H^{1,2} (\mathbb{B}^d ; \mathbb{S}^D)$,
       we say that the mapping $u \, \in \, $
        $H^{1,2} (\mathbb{B}^d ; \mathbb{S}^D)$
         is a weakly harmonic mapping (WHM) from $ \mathbb{B}^d$ into $\mathbb{S}^D$ 
          if the mapping $u$ satisfies the following P.D.E in the weak sense:
%
%\allowdisplaybreaks\begin{alignat}{2}
\begin{equation}
\left\{
\begin{array}{cll}
%1
- \triangle u 
 & \; = \; | \nabla u|^2 u 
  & \; \mathrm{in} \quad\mathbb{B}^d, 
\\
%2
u & \, = \, u_0 
 & \; \mathrm{on} \quad\partial \mathbb{B}^d.
\end{array}
\right.
\end{equation}
%\end{alignat}
%
Here \lq \lq P.D.E in the weak sense\rq\rq \enspace indicates that
 for $u \, - \, u_0 \, \in \,$ $\tc (\mathbb{B}^d ; \mathbb{R}^{D+1})$
  $=$ $\overline{C_0^\infty ( \mathbb{B}^d ; \mathbb{R}^{D+1})}^{H^{1,2}}$
   and any $\phi $ $\in$ $C_0^\infty (\mathbb{B}^d ; \mathbb{R}^{D+1})$,
\begin{align}
\allowdisplaybreaks
%1
\lint_{\mathbb{B}^d} & 
 \la \nabla u, \nabla \phi \ra \, dx 
  \, = \,
   \lint_{\mathbb{B}^d} |\nabla u|^2 \la u, \phi \ra \, dx
\label{EQ:Harm}
\end{align}
holds.
%
%\par
We employ the notation:
\allowdisplaybreaks
\begin{align}
%1
& \nabla u \, = \,
\left(\frac {\partial u^i}{\partial x_\alpha}\right)
 \quad (\alpha = 1, \ldots ,d ; i = 1, \ldots ,D+1), \quad
  \la u,v \ra \, = \, 
   \sum_{i=1}^{D+1} u^i v^i, 
\notag
\\
%2
& \la \nabla u, \nabla v \ra \, = \, 
 \sum_{\alpha=1}^d \sum_{i=1}^{D+1} 
  \frac {\partial u^i}{\partial x_\alpha} 
   \frac {\partial v^i}{\partial x_\alpha}, \quad
    | \nabla u | \, = \, \sqrt{\la \nabla u, \nabla u \ra}.
\end{align}
\par
It is easy to check that \eqref{EQ:Harm}
 is the Euler-Lagrange equations for the
  following variational problem of minimizing the Dirichlet energy
\begin{equation*}
 I[w] \, = \, \lint_{\mathbb{B}^d} | \nabla w |^2 \, dx
\label{EQ:Energy}
\end{equation*}
among mappings $w$ belonging to the admissible function class
\begin{equation}
H_{u_0}^{1,2} (\mathbb{B}^d ; \mathbb{S}^D) 
 \, := \, 
  \{w \, \in \, H^{1,2} (\mathbb{B}^d ; \mathbb{S}^D)
   ; \, w - u_0 
    \, \in \, \tc (\mathbb{B}^d ; \mathbb{R}^{D+1})\}.
\end{equation}
A work by R.Schoen and K.Uhlenbeck \cite{schoen-uhlenbeck-82}
 can read that any minimizing mapping of $I$
  in $H_{u_0}^{1,2} (\mathbb{B}^d ; \mathbb{S}^D)$
   denoted by $u_{\mathrm{min}}$ are smooth
    except possibly a closed set having at most $(d-3)$-Hausdorff dimension.
\par
An alternative approach to find a harmonic map,
 \lq \lq the heat flow method\rq \rq $\,$
  introduced by J.Eells and J.H.Sampson \cite{eells-sampson} is now a standard one.
They constructed a global harmonic heat flow 
 from any compact  manifold to any compact Riemannian manifold
  with non-positive sectional curvature.

Yielding to the approach, 
 when we let a positive number $T$ so large and
  a parabolic cylinder $Q(T)$ be $(0,T)$ $\times$ $\mathbb{B}^d$,
   we consider the heat flow:
\begin{equation}
\left\{
\begin{array}{rll}
%1
\dfrac {\partial u}{\partial t} & \, = \,
 \triangle u \, + \, | \nabla u |^2 u
  \quad & \mathrm{in} \;\; Q(T),
%\label{EQ:HHF}
\\[5pt]
%2
u(0,x) & \, = \, u_0(x)
 \qquad & \mathrm{at} \quad \{0\} \times \mathbb{B}^d,
\\[3pt]
%3
u(t,x) & \, = \, u_0(x)
 \qquad & \mathrm{on} \quad [0,T) \times \partial \mathbb{B}^d.
\end{array}
\right.
\label{EQ:HHF}
\end{equation}
We hereafter denote the time-slice mapping $v(t)$
 on $\mathbb{B}^d$ at a time $t$ of a mapping $v(t,x)$
  on $(0,T) \, \times \, \mathbb{B}^d$
   by $v(t)(x)$ $=$ $v(t,x)$.

Define three function spaces:
\allowdisplaybreaks\begin{align}
%1
L^\infty & \bigl(0,T;H^{1,2}(\mathbb{B}^d; \mathbb{S}^D)\bigr) 
\notag
\\
%2
&
\, := \, 
 \{ u | u \; \text{is measurable from} \;
  [0,T] \;\mathrm{to}\; H^{1,2} (\mathbb{B}^d ;\mathbb{S}^D) 
\notag
\\
%3
&
\quad \mathrm{with} \quad
 \underset {t \in (0,T)}{\esssup}
  ||u(t)||_{H^{1,2}(\mathbb{B}^d)} \, < \, + \infty \, \}, \hfill 
\notag
\\
%4
H^{1,2} & \bigl(0,T; L^2(\mathbb{B}^d;\mathbb{R}^{D+1}) \bigr) 
\notag
\\
%5
&
\; := \;
 \Bigl\{ u \, | \, u \; \mathrm{and} \;
  {\partial u}/{\partial t}
   \; \left(
    \text{a weak derivative of} \; u \right)
\notag
\\
%6
&
\quad \text{are second integrable from} 
 \: [0,T] \;\text{to}\; L^2 (\mathbb{B}^d;\mathbb{R}^{D+1}) 
  \Bigr\},
\notag
\\
%8
&
V (Q(T);\mathbb{S}^D) \; = \; 
 L^\infty \bigl(0,T;H^{1,2}(\mathbb{B}^d;\mathbb{S}^D)\bigr) \cap 
  H^{1,2} \bigl(0,T;L^2(\mathbb{B}^d;\mathbb{R}^{D+1})\bigr).
\notag
\end{align}
By virtue of topological obstruction, we generally have no hope
 about the existence of the classical solutions of the systems above 
  \eqref{EQ:HHF}.
So we need a weak formulation of it.
 For any given mapping $u_0 \, \in \,$ 
  $H^{1,2} (\mathbb{B}^d;\mathbb{S}^D)$,
   we call a mapping $u \, \in \, V(Q(T);\mathbb{S}^D)$
    a weakly harmonic heat flow {(\textit{WHHF})} provided
for any $\phi \in C_0^\infty ({Q(T)};\mathbb{R}^{D+1})$
\allowdisplaybreaks\begin{alignat}{2}
%1
&
\lint_{Q(T)} \Bigl( \la \frac {\partial u}{\partial t}, \phi \ra
 \, + \, \la \nabla u, \nabla \phi \ra
  \, - \, \Bigr. && \Bigl. %\kern-7mm
   \la u, \phi \ra |\nabla u|^2
    \Bigr) \, dz \, = \,0,
\label{EQ:1}
\\
%2
& 
u(t) \, - \, u_0 \in
 \tc (\mathbb{B}^d ;\mathbb{R}^{D+1}) 
  &&\quad \text{for almost every} \;
   t \, \in \, (0,T),
\label{EQ:2}
\\
%3
& \underset {t \searrow 0}{\lim}
 u(t) \, = \, u_0
  \quad && \quad \mathrm{in} \quad L^2 (\mathbb{B}^d ;\mathbb{R}^{D+1}).
\label{EQ:3}
\end{alignat}
Notice:
\begin{Thm}\label{THM:Wedge}
 \eqref{EQ:1} is equivalent to
\begin{equation}
\left\{
\begin{array}{ccl}
%1
& \dfrac {\partial u}{\partial t} \wedge u
 \, - \, \triangle u \wedge u \, = \, 0
  \quad & {\rm in}
   \quad \bigl(C_0^\infty ({Q(T)};\mathbb{R}^{D(D+1)/2}) \bigr)^*,
\label{EQ:Wedge}
\\
%2
& | u | \, = \, 1
 \quad & {\rm in} \quad {\rm a.e} \; z \in {Q(T)}.
\end{array}
\right.
\end{equation}
\end{Thm}
\par
To construct some WHHF, Y.Chen \cite{chen-89} plied
 the following penalty scheme
\allowdisplaybreaks\begin{align}
%1
\frac {\partial u_\lambda}{\partial t} & \, - \,
 \triangle u_\lambda \, + \, 
  \lambda ( |u_\lambda|^2 \, - \, 1) u_\lambda
   \; = \; 0
\label{EQ:Penalty}
\end{align}
and send $\lambda \nearrow \infty$.
Y.Chen and M.Struwe \cite{chen-struwe} has established 
 an existence and a partial regularity on a weakly harmonic heat flow
  between a compact Riemannian manifolds, and some later
   Y.Chen and F.H.Lin \cite{chen-lin} generalized it the case from
    a compact Riemannian manifold with a boundary to a compact Riemannian manifold.
Invoking a penalty approximation scheme, 
 they proved that it is smooth 
  except a set called \lq\lq singular set\rq\rq 
   ~having at most the finite $d$-dimensional Hausdorff measure
    with respect to the parabolic metric.
X.Cheng \cite{cheng} showed that the time slice of singular set
 has at most $(d-2)$-dimensional Hausdorff measure at {\it{every}} time moment, 
instead of {\it{almost}} every time.
Inspired by the work of L.Caffarelli, R.Kohn and L.Nirenberg
 \cite{caffarelli-kohn-nirenberg} or L.C.Evans \cite{evans-91}
  or F.H\'{e}lein \cite{helein},
   Y.Chen, J.Li and F.H.Lin \cite{chen-li-lin} and 
    M.Feldman \cite{feldman} discussed
     a partial regularity for a WHHF in a certain function class.
Intrinsically, they assumed that their harmonic heat flow
 has a monotonicity for the scaled energy and a local energy
  estimate.
Then they showed that such a flow possibly 
 has a singular set with {\it{zero}} $d$-dimensional
  Hausdorff measure with respect to the parabolic metric.
\par
As an improved approach to construct a WHHF,
 I will propose a new approximate evolutional scheme
  said to be the Ginzburg-Landau heat flow  and abbreviated by GLHF.
To explain Ginzburg-Landau heat flow,
 we introduce smooth functions $\chi (t)$ 
  and $\kappa (t)$ by
\allowdisplaybreaks\begin{align}
%1
&
\chi (t) \; = \;
\begin{cases}
%1
t & \quad (t \, < \, 2)
\\
%2
3 & \quad (t \, \ge \, 4),
\end{cases}
\quad \chi \, \le \, 3,
\notag
\\
%2
&
\kappa (t) \; = \; \arctan (t)/ \pi.
\end{align}
For a mapping $u_0 \, \in \, H^{1,2} ( \mathbb{B}^d ; \mathbb{S}^D )$
 $\cap$ $H^{2,2} ( \mathbb{B}^d\setminus {B}_{1-\delta_0} (0) ; \mathbb{S}^D )$
with a positive number sufficiently small $\delta_0$,
 a GLHF is designated by solution of the systems:
%
%\begin{equation}
%\left\{
\begin{eqnarray}%{rlc}
%1
\dfrac{\partial u_\lambda}{\partial t} \, - \,
 \triangle u_\lambda
  \, + \, \lambda^{1-\kappa}
   \dot \chi \bigl( \bigl( | u_\lambda |^2 \, - \, 1 \bigr)^2 \bigr)
    \bigl( | u_\lambda |^2 \, - \, 1 \bigr) u_\lambda 
     \; = \; 0
      \quad \mathrm{in} \quad Q(T),
\label{EQ:GLHF}
\\[2mm]
%2
u_\lambda
 \; = \; u_0
  \quad \mathrm{on} \quad \partial Q(T).
%\end{array}
%\right.
\label{EQ:GLHF-Bdry}
\end{eqnarray}
%\end{equation}
%
If you notice that the non-linear term of 
 $\lambda^{1-\kappa}$
  $\dot \chi \bigl( \bigl( | u_\lambda |^2 \, - \, 1 \bigr)^2 \bigr)$
   $\bigl( | u_\lambda |^2 \, - \, 1 \bigr) u_\lambda$
    is bounded, 
     Banach's fixed point theorem can state
      the unique existence of the mapping $u_\lambda$ on $Q(T)$ with
\begin{enumerate}
%1
\item[(a)]
$u_\lambda \, \in \, C^\infty (Q(T))$,
%2
\item[(b)]
\text{\eqref{EQ:GLHF} is fulfilled in } $Q(T)$,
%3
\item[(c)]
$u_\lambda (t,x) \, - \, u_0 (x)$
 $\in$
  $\tc (\mathbb{B}^d;\mathbb{R}^{D+1})$ 
$\;\, \text{for almost every}$ 
 $\; t \;\mathrm{in}\; (0,T)$,
%4
\item[(d)]
$
\lim_{t \searrow 0}
 || u_\lambda (t, \circ) \, - \, u_0 (\circ) ||_{L^2 (\mathbb{B}^d)}
  \; = \; 0$.
\end{enumerate}
\label{DEF:GLHF}
\par
We mention it in Theorem \ref{THM:Existence-GLHF} in
 p.p \pageref{P:GLHF}.
I emphasis that the time dependent parameter $\lambda^{1-\kappa}$
 in \eqref{EQ:GLHF} easily leads to 
  $\int_{Q(T)} \lambda^{1-\kappa}$
   $( | u_\lambda |^2 \, - \, 1 )^2$ $dz$ $ \le $
    $C/\log \lambda$; This makes us handle with the nonlinear term.
     The reader should inquire Theorem \ref{THM:Energy-Estimate}.
\par
We give a few remarks about the hypothesis on the mapping $u_0$ above.
 First of all, on account of a topological obstruction,
  a class of mappings $C^1 (\overline{\mathbb{B}}^d ; \mathbb{S}^{d-1})$
   is empty if the degree of their restriction to the boundary 
    doesn't vanish, 
whereas so isn't $H^{1,2} ({\mathbb{B}}^d ; \mathbb{S}^{d-1})$
 if $d$ is more than or equal to $3$.
For instance, 
 for any mapping $\phi$ $\in$ $C^1$ $(\partial{\mathbb{B}}^d ;
  \mathbb{S}^{d-1})$,
   the mapping $\phi (x/{|x|})$ belongs to 
    $H^{1,2} ({\mathbb{B}}^d ; \mathbb{S}^{d-1})$    
as long as $d$ is more than or equal to $3$.
 F.Betuel and X. Zheng \cite{bethuel-zheng}
  systematically studied a density result
   of various Sobolev mappings between two Riemannian manifolds.
They pointed out that any map in $H^{1,2} ({\mathbb{B}}^d ; \mathbb{S}^{d-1})$
 is approximated by the   mapping that is smooth except finite points.
Second the further imposition of the initial and boundary mapping:
 the mapping $u_0$ belongs to $H^{2,2}$ near the boundary,
  is necessary for the legitimacy of 
   Theorem \ref{THM:Energy-Estimate} (Energy estimates),
    Theorem \ref{THM:EDE} (Energy decay estimates),
     Corollary \ref{COR:PPI} (Monotonicity for the scaled energy) 
      and Theorem \ref{THM:Mon} (Parabolic Pokhojaev inequality).
\par
The aim of the paper is to construct a WHHF with
\renewcommand{\theenumi}{\roman{enumi}}
\begin{enumerate}
\item
a global energy inequality,
\item
a monotonicity for the scaled energy,
\item
a reverse Poincar\'e inequality.
\end{enumerate}
We should come to our mind that (\romannumeral 2) and (\romannumeral 3) are indispensable tools 
 in a regularity or a partial regularity theory on the
  solutions of various elliptic and parabolic equations.
This paper first establishes a crude bound,
 a maximal principle, a global energy inequality,
  a monotonicity inequality for the scaled energy 
   and the reverse Poincar\'e inequality for GLHF; 
Next we show that the GLHF converges to a WHHF as $\lambda \nearrow \infty$
  (modulo a subsequence of $\lambda$) and the flow also satisfies
   (\romannumeral 1), (\romannumeral 2) and (\romannumeral 3).
Thereafter by making the best of a monotonicity inequality for the scaled energy
 and the reverse Poincar\'e inequality,
  we will prove that 
   the WHHF actually is smooth except on a small set
    called \lq\lq \bf{singular set}.\rq\rq \rm \enspace
     More precisely we assert
%%%%%%%%%%%%%%%%%%%%%%%%%%%%%%%%%%%%%%%%%%%
\begin{Thm}{\rm{(Partial Regularity).}}\label{THM:Main-1}
Let $d$ be a positive integer larger than or equal to $3$.
For a mapping $u_0 \, \in \, H^{1,2} ( \mathbb{B}^d ; \mathbb{S}^D )$
 $\cap$ $H^{2,2} ( \mathbb{B}^d\setminus {B}_{1-\delta_0} (0) ;
  \mathbb{S}^D )$
   with a sufficiently small positive number $\delta_0$,
there exists a WHHF
 and it is smooth on a certain open set in $Q (T)$ whose 
  compliment has the finite $(d-\epsilon_0)$-dimensional Hausdorff measure
   with respect to the parabolic metric,
where $\epsilon_0$ is a small positive number depending only on $u_0$ and $d$.
The WHHF constructed above also holds
\allowdisplaybreaks\begin{align}
%1
&
{\rm{(\romannumeral1)}} \hfill \quad
 \lint_{T_1}^{T_2} \, dt \lint_{\mathbb{B}^d}
  \left| \frac{\partial u}{\partial t} \right|^2 \, dx
\, + \, \frac 12 \lint_{\mathbb{B}^d} |\nabla u (T_2,\cdot)|^2 dx
\; \le \; \frac 12
 \lint_{\mathbb{B}^d} |\nabla u (T_1,\cdot)|^2 dx
\label{INEQ:Main-Energy-Estimate-1}
\\
\intertext{for almost every time $T_1$ and $T_2$ with
 $0$ $\le$ $T_1$ $\le$ $T_2$ $\le$ $T$,}
%2
&
{\rm{(\romannumeral2)}} \hfill \quad
 \lint_{R_1}^{R_2} \frac {d\rho}{\rho^{d-1}}
  \lint_{t_0 - (2\rho)^2}^{t_0 - \rho^2} \, dt \lint_{\mathbb{B}^d }
 \left| \frac {\partial u}{\partial t}
  \; + \;
   \frac {x \, - \, x_0}{2(t \, - \, t_0)} \cdot \nabla u
		    \right|^2
    \exp \Bigl( \frac{|x \, - \, x_0|^2}{4(t \, - \, t_0)} \Bigr) \, dx
\notag\\
%3
&
\qquad \, + \, \frac 1 {2{R_1}^d}
  \lint_{t_0 - (2{R_1})^2}^{t_0 - {R_1}^2} \, dt
%   \kern-3mm \, dt
    \lint_{\mathbb{B}^d} | \nabla u |^2 
     \exp \Bigl( \frac{|x \, - \, x_0|^2}{4(t \, - \, t_0)} \Bigr) \, dx
\label{INEQ:Main-Energy-Estimate-2}
\\
%4
&
\qquad \; \le \;
 \frac 1 {2{R_2}^d}
  \lint_{t_0 - (2{R_2})^2}^{t_0 - {R_2}^2} \, dt
%   \kern-3mm \, dt
    \lint_{\mathbb{B}^d} | \nabla u |^2
     \exp \Bigl( \frac{|x \, - \, x_0|^2}{4(t \, - \, t_0)} \Bigr) \, dx
\, + \, \frac {C_{\mathsf{M}}}2 (R_2^2 - R_1^2)
\notag
\end{align}
for any positive numbers $R_1$ and $R_2$ with $R_1 \le R_2$ 
 and an arbitrary point $z_0$ $=$ $(t_0,x_0)$ in $Q(T)$
  satisfying $t_0 - (2R_2)^2$ $> 0$,
   where $C_{\mathsf{M}}$ is a certain positive
    constant depending only on $d$, $z_0$ and $u_0$,
and in addition
\allowdisplaybreaks\begin{align}
%1
&
{\rm{(\romannumeral3)}} \hfill \quad
 \lint_{P_R (z_0)} | \nabla u (z) |^2 \, dz
  \; \le \; \frac {C}{R^2} \lint_{P_{2R} (z_0)}
   | u (z) \, - \, a (t) |^2 \, dz
\end{align}
for any $t$-variable second integrable mapping $a$ $=$ $a(t)$ 
 $=$ $( a^i (t))$ $( i = 1,2,\ldots, D+1)$
 and any parabolic cylinder 
  $P_{2R}$ $(z_0)$ compactly contained in $Q(T)$.
\end{Thm}
%%%%%%%%%%%%%%%%%%%%%%%%%%%%%%%%%%%%%%%%%%%KOKO
\par
We next discuss an asymptotic behaviour of the WHHF as $t \nearrow \infty$;
 The final part of the paper demonstrates
  that our WHHF with a constant boundary value
   converges strongly to a constant as $t \nearrow \infty$
    in $H^{1,2} (\mathbb{B}^d ; \mathbb{S}^D)$
     as long as the dimension $d$ is more than or equal to $3$. 
A surprising example by K-C.~Chang and W-Y.~Ding \cite{chang-ding}
 tells us that the same result doesn't hold in $d=2$.
  They constructed a smooth harmonic heat flow $u$ $=$ $u(t,x)$
   under a smooth initial condition $u_0$ with $u_0 (\mathbb{B}^2)$
    $=$ $\mathbb{S}^2$ 
     and a boundary condition $\left. u_0 \right\vert_{\partial \mathbb{B}^2}$
      $=$ $\mathit{a \; constant}$, and 
       showed that it \textit{does not} converge to any WHM
        in $C^0(\mathbb{B}^2)$ as $t \nearrow \infty$.
%
%%%%%%%%%%
%
\par
Their result entails that $u (t)(\mathbb{B}^2)$ $=$ $\mathbb{S}^2$ holds
 at each time $t$ in $(0,\infty)$.
On the other hand, non-existence result
 by L.Lemaire \cite{lemaire-1} and \cite{lemaire-2} stands for
\begin{Thm}\label{THM:Const}
Let a mapping $\phi$ $:$ $\mathbb{B}^2 \, \longrightarrow \, \mathbb{S}^2$
 be harmonic with $\left. \phi \right\vert_{\partial \mathbb{B}^2}$
  $=$ $\mathit{a \; constant}$. 
   Then the mapping $\phi$ must be the same constant.
\end{Thm}
Thus the harmonic heat flow constructed above is not homotopic to $\phi$
 nevertheless it has the same boundary condition.
  So we can not expect that the harmonic heat flow $u$ converges to \it{a constant} \rm
   as $t \nearrow \infty$ in $C^0 (\mathbb{B}^2)$.
In a nut shell the square sum of the gradients of the harmonic heat flow 
 by K-C.~Chang and W-Y.~Ding,
  converges to the Dirac measure as $t \nearrow \infty$.
   When $d=2$, we should notice that the WHHF is unique
    due to A. \cite{freire} under the following criterion:
\begin{equation}
t \; \longmapsto \; \lint_{\mathbb{B}^d} | \nabla u(t) |^2 \, dx
\label{PROP:NI}
\end{equation}
is non-increasing in $[0,\infty)$.
 On the other hand, we know the existence of infinite many WHHFs
  not satisfying the criterion above.
   This is by M.Bertsch, R.Dal Passo and R.~Vd Hout
    \cite{bertsch-dalpasso-Rvdhout}.
\par
When $d$ is greater than or equal to $3$, the equator map 
 ${x}/{|x|}$ $\in$ $H^{1,2} (\mathbb{B}^d ;
 \mathbb{S}^{d-1})$\footnote{This map is also WHM and
  even the absolute minimizer. We refer it to H.~Brezis, J.M.Coron and
   E.H.Lieb \cite{brezis-coron-lieb}}
     reveal us the WHHF possessing a singularity.
      Thus any WHHF may break any given topology.
In $d=3$, the uniqueness of WHHF fails even if \eqref{PROP:NI}
 does hold and $u_0$ is a WHM by J.M.Coron \cite{coron}
  or any stationary by M.-C.Hong \cite{hong}
   or a minimizer by M.Bertsch, R. Dal Passo and A. Pisante
    \cite{bertsch-dalpasso-pissante}.
\par
Our alternative main theorem is
\begin{Thm}\label{THM:Main-2}
Let a mapping $u$ be the WHHF obtained by the limit of GLHF as 
 $\lambda \nearrow \infty$ {\rm{(}} modulo of subsequence of $\lambda$ {\rm{).}}
Then we can extend the WHHF from $Q(T)$ to $Q(\infty)$
 and moreover 
  the mapping $u(t)$ converges strongly to a constant in $H^{1,2} (\mathbb{B}^d)$
   as $t \nearrow \infty$ if $u (t) \vert_{\partial \mathbb{B}^d}$
    $=$ \textit{the \,constant} \textrm at almost all $t>0$.
\end{Thm}
This result can be regarded as a parabolic analogue of non-existence result
 of harmonic mappings.
The proofs of Theorem \ref{THM:Main-1} and Theorem \ref{THM:Main-2} will be given
 in the final part of the paper.
\vskip 9pt
\par
We close this introduction by enumerating a glossary of notation:
\vskip 18pt
\vskip 9pt
\begin{centerline}
{Notation}
\end{centerline}
\vskip 9pt
\renewcommand{\labelenumi}{(\roman{enumi})}
\begin{enumerate}
\item
 $\mathbb{B}^d$
  $\, = \,$ 
   $\{x = (x_1,x_2,\ldots,x_d) \in \, \mathbb{R}^d;
    |x| = \sqrt{\sum_{\alpha=1}^d (x_\alpha)^2 }\, < 1\, \}$
%\item

and $\; \partial \mathbb{B}^d$ is the boundary of $\mathbb{B}^d.$
\item
${C}^d$
 $\, = \,$
  $\{x = (x_1,x_2,\ldots,x_d) \in \, \mathbb{R}^d; \;
   |x_\alpha| \, < \, 1 \;
    ( \alpha \, = \, 1,\ldots,d)
     \}$.
\item
 ${Q(T)}=(0,T) \times \mathbb{B}^d$.
%\item
 $\; \partial {Q(T)} \, = \,$
  $[0,T) \times \partial \mathbb{B}^d$
   $\cup$
    $\{0\} \times \mathbb{B}^d$.
\item
$\mathbb{S}^D  \, = \, $ 
 $\{y = (y^1,y^2,\ldots, y^{D+1}) \, \in \, \mathbb{R}^{D+1};$ 
  $|y| \, = \, \sqrt{\sum_{i=1}^{D+1} (y^i)^2} \, = \, 1\}.$
\item
For any points $z_1 = (t_1,x_1)$ and $z_2 = (t_2,x_2)$,
 the parabolic metric $d(z_1,z_2)$ means
  $|t_1 \, - \, t_2|^{1/2}$ $\, + \,$ $|x_1 \, - \, x_2|$.
\item
For any positive integer $n$, $Q_n$ $=$
 $[1/n,T-1/n] \times \overline{B_{1-1/n} (0)}$
  and $d_n$ $=$ ${d}$ $(Q_n, \partial Q(T))$.
\item
For a Lebesgue measurable subset $A$ in $\mathbb{R}^d$ or ${Q(T)}$,
 $|A|$ denotes the $d$ or $(d+2)$-dimensional Lebesgue measure
  of $A$.
%%
%\item
%$[*]$ is the Gauss symbol.
%
\item
Set points $x$ $=$ $(x_1,x_2,\ldots,x_d)$,
 $x_0$ $=$ $(x_{0,1},x_{0,2},\ldots,x_{0,d})$ and
  $z_0 \, = \, ( t_0 , x_0) \in \mathbb{R}^{d+1}$;
   Then indicate a ball $B_r (x_0)$, a cube $C_r (x_0)$,
    a parabolic cylinder $P_r (z_0)$ 
     and a  parabolic cube $D_r (z_0)$ by
\allowdisplaybreaks\begin{align*}
%1
B_r (x_0) &
 \, = \, \{ x \in \mathbb{R}^d \, ; \,  | x \, - \, x_0 | < r \}, \\
%2
C_r (x_0) &
 \, = \, \{ x \in \mathbb{R}^d \,; \, | x_i \, - \, x_{0,i} | < r,
            \; ( i = 1,2,\ldots,d) \},
\\
%3
P_r (z_0) & \, = \, 
 \{ z=(t,x) \in Q(T)  \, ; \, t_0 - r^2 < t < t_0 + r^2,
  | x \, - \, x_0 | < r\},
\\
%4
D_r (z_0) & \, = \,
 (t_0-r^2,t_0+r^2) \times C_r (x_0).
\end{align*}
In $B_r (x_0)$, $C_{r} (x_0)$, $P_r (z_0)$ and $D_r (z_0)$,
  the points of $x_0$ and $z_0$ will be often abbreviated
   when no confusion may arise.
\item
For a set $A \, \subset \mathbb{R}^d$ and
 $0 \le d^\prime < \infty$,
  we define the $d^\prime$-dimensional Hausdorff measure
\footnote{This definition is not it on
          the usual Hausdorff measure. However it is sufficient
           to our end.}
with respect to the parabolic metric by
\allowdisplaybreaks\begin{align*}
%1
&
\mathcal{H}^{(d^\prime)} (A) 
 \; = \; \lim_{R \searrow 0} \mathcal{H}_{R}^{(d^\prime)} (A)
\\
%2
&
\mathrm{with} \quad
 \mathcal{H}_{R}^{(d^\prime)} (A) \; := \;
  \inf_{\mathrm{covering}}
   \Bigl\{ \sum_i R_i^{d^\prime} ;\; 
    A \subset \bigcup_i P_{R_i} (z_i), \, 0 < R_i < R \Bigr\}.
\end{align*}
\item
Letter $C$ denotes a generic constant.
 By the letter $C(B)$, it means that a constant depends only on a parameter $B$.
\item
$\kappa (t)$ is $\arctan (t) / \pi$.
\item
$[*]$ is the Gauss symbol.      
\end{enumerate}
\vskip 9pt
\begin{centerline}
{Functions and Derivatives}
\end{centerline}
\vskip 9pt
\begin{enumerate}
\item
For vectors $u$ and $v$$\in \mathbb{R}^{D+1}$,
 $\la u,v \ra$
  $\, = \, $
   $\sum_{i=1}^{D+1} u^iv^i$
    and
     $|u|$ $\, = \,$ $\la u,u \ra^{1/2}$.
\item
For a map $v: \mathbb{R}^d \to \mathbb{R}^{D+1}$
 the gradient matrix of $v$ is defined by 
  $\nabla v$ $\, = \, $$\left( \partial v^i / \partial x_\alpha \right)$
   $(\alpha =1, \ldots ,d \, ; \, i=1,\ldots ,D+1)$,
    while $\partial_t v$ is by ${\partial v}/\partial t$.
In addition, for $x$ $=$ $(x_1, x_2, \ldots, x_d)$,
 $x\!\cdot\! \nabla$ and $\nabla_\nu$ respectively denotes 
  $\sum_{\alpha=1}^d$ $x_\alpha \, {\partial}/{\partial x_\alpha}$
   and $\nu \cdot \nabla$, 
    a symbol $\nu$ denoting the outward normal unit vector on the boundary of a discussed domain.
\par
On the contrary, differentials
 $\nabla_\tau$ $=$ $( \nabla_{\tau_i} )$ $( i = 1,2,\ldots, d-1 )$
  indicates by $\nabla$ $-$ $\nu$ $\nu \cdot \nabla$ 
   and $\triangle_\tau$ $=$ $\sum_{j=1}^{d-1}$ $\nabla_{\tau_i}^2$.
\par
For any positive integer $k$,
 a higher differentials $\nabla^k$ and $\nabla_\tau^k$ mean
\begin{align}\allowdisplaybreaks
%1
\nabla^k & \; = \;
\sum
_{\begin{subarray}{l}
k_1 \, + \, \cdots \, + \, k_{d-1} \; = \; k \\
0 \, \le \, k_1, \ldots, k_{d-1} \, \le \, k
\end{subarray}}
\nabla_{x_1}^{k_1} \, \cdots \, \nabla_{x_{d-1}}^{k_{d-1}},
\notag
\\
%2
\nabla_\tau^k & \; = \;
\sum
_{\begin{subarray}{l}
k_1 \, + \, \cdots \, + \, k_{d-1} \; = \; k \\
0 \, \le \, k_1, \ldots, k_{d-1} \, \le \, k
\end{subarray}}
\nabla_{\tau_1}^{k_1} \, \cdots \,\nabla_{\tau_{d-1}}^{k_{d-1}}
\notag
\end{align}
\item
For the gradient matrices of $u$ and $v$,
 namely $\nabla u$ and $\nabla v$,
  we mean $\la \nabla u, \nabla v \ra$
   $\, = \, $
    $\sum_{\alpha=1}^d \sum_{i=1}^{D+1}$
     $\nabla_\alpha u^i \nabla_\alpha v^i$
      and
       $|\nabla u|^{1/2}$
        $\, = \,$$\la \nabla u, \nabla u \ra^{1/2}$.
\item
For GLHF $u_\lambda$, 
 we call $\mathbf{e}_\lambda$ 
  the Ginzburg-Landau energy density 
   given by $|\nabla u_\lambda|^2/2 \, + \, \lambda^{1-\kappa} 
    \chi ((| u_\lambda |^2 \, - \, 1 )^2)/4$.
\item
For a mapping $u(t,x)$ on $(t,x)$ $\in$ $Q(T)$,
 a time-slice mapping $u(t)$ at a time $t$, denotes
  $u(t)(x)$ $=$ $u(t,x)$.
\end{enumerate}
\vskip 9pt
\begin{centerline}
{Function spaces}
\end{centerline}
\vskip 9pt
%Let $p$ be a number more   than or equal to $1$.
\begin{enumerate}
\item
What a function $f$ belongs to 
 $C^0 (\mathbb{B}^d)$
  or $C^0 (Q(T))$ is
   that the function $f$ is continuous on 
    $\mathbb{B}^d$ or $Q(T)$.
\item
$C_0^\infty (\mathbb{B}^d)$
 or $C_0^\infty ({Q(T)})$ is respectively
  the space of infinite differentiable function 
   with a compact support in $\mathbb{B}^d$ or ${Q(T)}$ and
    $(C_0^\infty)^*$ is the dual space of it.
\item
We say that a function $f$ belongs to H\"older space on $\overline{Q(T)}$
 if there is a positive constant $C$ and a positive number
  $\alpha_0$ $(0<\alpha_0<1)$ such that
$$
|f(z_1) \, - \, f(z_2) \, | \le \, C d(z_1,z_2)^{\alpha_0}
 \quad \mathrm{for \enspace any} \, \, z_1, z_2 \, \in \, 
  \overline{Q(T)}.
$$
The semi-norm of such a function $f$ is given by \enspace
\begin{equation*}
[f]_{\mathrm{C}^{\alpha_0} (\overline{Q(T)})} =
 \sup_{
  \begin{subarray}{lc}
   z_1,z_2 \in \overline{Q(T)}, \\
   z_1 \ne z_2
  \end{subarray}}
\frac {|f(z_1) \, - \, f(z_2)|}{d(z_1,z_2)^{\alpha_0}}.
\end{equation*}
\item
$C^{2,\alpha_0} (\overline{Q(T)})$
 $\; = \;$
  $\{ f \, \in \, C^0 (\overline{Q(T)}) \; ; \; \nabla f,
   \nabla^2 f \, \mathrm{and} \,
    \partial_t f \;\, 
     \text{is \thinspace continuous \thinspace on $\overline{Q(T)}$} \;
      \}$
of which the norm 
\begin{equation*}
\begin{split}
%1
||f||_{\mathrm{C}^{2,\alpha_0} (\overline{Q(T)})} 
 \; = \; &
  \sup_{z \in (\overline{Q(T)})} | f (z) |
   \, + \, 
    \sup_{z \in (\overline{Q(T)})} | \nabla f (z) |
\\
%2
&
\, + \,
 \sup_{z \in (\overline{Q(T)})} | \nabla^2 f (z) |
  \, + \,
   \sup_{z \in (\overline{Q(T)})} | \partial f / \partial t (z) |
\\
%3
&
\, + \, [\nabla^2 f]_{\mathrm{C}^{\alpha_0} (\overline{Q(T)})}
 \, + \, [\partial f / \partial t]_{\mathrm{C}^{\alpha_0/2} (\overline{Q(T)})}
\end{split}
\end{equation*}
is finite.
\item
$L^p (\mathbb{B}^d)$ or $L^p ({Q(T)})$ respectively means the space of 
 the $p$th summable function 
  on $\mathbb{B}^d$ or ${Q(T)}$ with the norm of
   $||f||_{L^p (\mathbb{B}^d)} \, = \, $
    $\bigl(\int_{\mathbb{B}^d} | f |^p dx \Bigr)^{1/p}$
     or
      $||f||_{L^p ({Q(T)})} \, = \, $
       $\Bigl(\int_{Q(T)} | f |^p dz \Bigr)^{1/p}$.
On the contrary, $L^\infty (\mathbb{B}^d (Q(T)))$
 is the space of any summable function so that
  the norm of $\Vert f \Vert_{L^\infty (\mathbb{B}^d)}$
   is $\sup_{x \in \mathbb{B}^d}$ $| f(x)|$ and
    the one of $\Vert f \Vert_{L^\infty (Q(T))}$ 
     $\sup_{z \in Q(T)}$ $| f(z)|$.

\item
$H^{1,2} (\mathbb{B}^d)$
 $\,=\,$ 
  $\{ f \in L^2( \mathbb{B}^d ) \, ; \, \partial f/ \partial x_\alpha \in L^2( \mathbb{B}^d )$,
   $( \alpha  = 1,\ldots ,d)\}$.
%
%The norm of $H^{1,2} (\mathbb{B}^d)$ is 
%%\newline
%  $||f||_{H^{1,2} (\mathbb{B}^d)} \, = \, $
%   $\Bigl( \int_{\mathbb{B}^d} | f |^2 dx \, 
%    + \, \int_{\mathbb{B}^d} | \nabla f |^2 dx \Bigr)^{1/2}$.
%
\item
$H^{1,2} ({Q(T)})$
 $\,=\,$ 
  $\{ f \in L^2( {Q(T)} ) \, ; \, \partial f/ \partial x_\alpha, \partial_t f \in L^2( {Q(T)} )$, 
   $( \alpha  = 1,\ldots ,d)\}$.
%
%The norm of $H^{1,2} ({Q(T)})$ is given by
%\newline
%  $||f||_{H^{1,2} ({Q(T)})} \, = \, $
%   $\Bigl( \int_{Q(T)} | f |^2 \, dz 
%     \, + \, \int_{Q(T)} |\nabla f |^2 \, dz 
%      \, + \, \int_{Q(T)} |\partial_t f |^2 dz \bigr)^{1/2}$.
%
\item
$H^{2,p} ({Q(T)})$
 $\,=\,$
  $\{ f \in L^p( {Q(T)} ) \, ; \, \partial f/ \partial x_\alpha, 
   \partial^2 f/ \partial x_\alpha \partial x_\beta,
    \partial_t f \in L^p( {Q(T)} )$,
     $( \alpha, \beta  = 1,\ldots ,d)\}$, 
where a number $p$ is  more than or equal to $1$.
\end{enumerate}
\par
In the following, 
 let $X(\mathbb{B}^d ( \mathrm{or} \; {Q(T)} ))$ be a Banach space
  on $\mathbb{B}^d ( \mathrm{or} \; {Q(T)})$.
\renewcommand{\labelenumi}{(\roman{enumii})}
\begin{enumerate}
\setcounter{enumii}{9}
\item
If a function $f$ belongs to 
 $X_{\mathrm{loc}} (\mathbb{B}^d ( \mathrm{or} \; {Q(T)} ))$,
  this means that the function $f$ is of
   $X_{\mathrm{loc} (\Omega)}$ for any set $\Omega$
    compactly contained in $\mathbb{B}^d ( \mathrm{or} \; {Q(T)} )$.
\item
$\overset \circ {X}(\mathbb{B}^d)$
 is the subspace of $X(\mathbb{B}^d)$
  whose element vanishes on $\partial \mathbb{B}^d$
   in the trace sense.
\setcounter{enumii}{10}
\item
$X \bigl(\mathbb{B}^d ( \mathrm{or} \; {Q(T)} ) ;\mathbb{R}^{D+1}
    \bigr)$
 $\,=\,$
  $\{ u = (u^i) : \mathbb{B}^d ( \mathrm{or} \; {Q(T)} ) \to
     \mathbb{R}^{D+1}$
   $; u^i \in \, X\bigl(\mathbb{B}^d (\mathrm{or} \; {Q(T)}) \bigr),$
    $\, (i=1, \ldots ,D+1)\}$.
\setcounter{enumii}{11}
\item
$X \bigl(\mathbb{B}^d (\mathrm{or} \; {Q(T)}) ; \mathbb{S}^D \bigr)$
 $\,=\,$
  $\{ u = (u^i) \, \in \, X(\mathbb{B}^d (\mathrm{or} \; {Q(T)}) ;
   \mathbb{R}^{D+1})$;
    $|u| \, =  \, 1 \, \, \mathrm{a.e} \, \, x \in \mathbb{B}^d$
     $( \mathrm{or \enspace a.e} \, \, z \in {Q(T)}),$
      $\, (i=1, \ldots ,D+1)\}$.
\setcounter{enumii}{12}
\item
If we say that a mapping 
 $u$ $=$ $u(t,x)$ on $(0,T) \times \mathbb{B}^d$ belongs to

  $Y (0,T ; X(\mathbb{B}^d ; \mathbb{R}^{D+1} (\mathrm{or} \; \mathbb{S}^D)))$, it means
   $u(t)$ $\in$ $X$ and $|| u(t) ||_X$ $\in$ $Y$,
where a symbol $|| \cdot ||_X$ is the equipped norm in a normed space $X$.
\setcounter{enumii}{13}
\item
Let a mapping $v$ be any mappings belonging to $H^{1,2}(\mathbb{B}^d ;\mathbb{S}^D)$.
 If a mapping $w \, \in \, $$H^{1,2} (\mathbb{B}^d ; \mathbb{S}^D)$
  and $w \, - \, v$
   $\, \in \, \tc$
    $(\mathbb{B}^d ; \mathbb{R}^{D+1})$,
     we then call the mapping $w$ belong to
      $H_v^{1,2} (\mathbb{B}^d ; \mathbb{S}^D)$.
\setcounter{enumii}{14}
\item
$V ({Q(T)};\mathbb{S}^D)$ $\; = \;$
  $L^\infty \bigl(0,T;H^{1,2}(\mathbb{B}^d;\mathbb{S}^D)\bigr) \cap$ 
   $H^{1,2} \bigl(0,T;L^2(\mathbb{B}^d;\mathbb{R}^{D+1})\bigr)$.
\end{enumerate}

%% file: glhf_hhf_arch.tex
%
%#! platex hhf
%
\setcounter{chapternumber}{2}\setcounter{equation}{0}
\renewcommand{\theequation}%
             {\thechapternumber.\arabic{equation}}
\section{\enspace GLHF.}
\label{P:GLHF}
The chapter is devoted to the study of GLHF.
 We prove that a GLHF satisfies 
  a crude bound, a maximal principle, a few energy inequalities,
   a monotonicity inequality for the scaled energy
    and finally a hybrid type inequality.
\subsection{\enspace Properties on GLHF.}
As alluded in Chapter 1, we briefly discuss how to construct the GLHF.
 We first assume that the initial and boundary mapping $u_0$
  is smooth in $\overline{Q(T)}$ in Theorem \ref{THM:Existence-GLHF},
   Lemma\ref{LEM:CB}, Theorem \ref{THM:Max}, Theorem \ref{THM:Energy-Estimate},
    Theorem \ref{THM:EDE}, Corollary \ref{COR:PPI} and
     Theorem \ref{THM:LEI}:     
      We state
\begin{Thm}{\rm{(The existence of GLHF)}.} \label{THM:Existence-GLHF}
Give a smooth mapping $u_0$ between $\overline{Q(T)}$ and $\mathbb{R}^{D+1}$.
 Then there exists the classical solution $u_\lambda$ $\in$ $C^\infty$
  $(\overline{Q(T)})$ to \eqref{EQ:GLHF} under $u_\lambda$ $=$
   $u_0$ on $\partial Q(T)$ such that when we set
    $u_\lambda^{(1)}$ and $u_\lambda^{(2)}$ respectively
     the classical solution to \eqref{EQ:GLHF}
      with the smooth boundary and initial condition 
       $u_0^{(1)}$ and $u_0^{(2)}$ on 
        $\partial Q(T)$,
we have
\allowdisplaybreaks\begin{align}
%1
&   
|| u_\lambda^{(1)} \, - \, u_\lambda^{(2)} ||_{H^{2,2} (Q(T)}
\label{INEQ:Schauder}
\\
%2
&
\; \le \;
 C(\lambda,T)
  ( 1 \, + \, 
   || u_0^{(1)} \, -  \, u_0^{(2)} ||_{H^{1,2} ( \mathbb{B}^d )}
    \, + \,
     || u_0^{(1)} \, -  \, u_0^{(2)} ||_{H^{2,2} ( \mathbb{B}^d\setminus {B}_{1-\delta_0} (0)})                                                                                 \
\notag
\end{align}
with a positive number $\delta_0$ sufficiently small.
\end{Thm}
\vskip 9pt
\noindent{\underbar{{Proof of Theorem \ref{THM:Existence-GLHF}}}.}
\rm\enspace
\vskip 6pt
\par
A routine work by means of 
 a Duhamel's formula and a contraction mapping theorem
  implies the unique classical solution to
   \eqref{EQ:GLHF} on $Q( t_\lambda )$
    for a small positive number $t_\lambda$ possibly depending on $\lambda$.
We repeat the argument above to extend our solution 
 to the time interval $[ t_\lambda, 2t_\lambda )$.
  Continuing after finite steps, we eventually come up with the classical solution
   to \eqref{EQ:GLHF} on $Q(T)$.
\par
Moreover, since the mapping $u_\lambda^{(1)} \, - \, u_\lambda^{(2)}$
 is the solution to
\begin{equation}
\left\{
\begin{array}{ll}
%1
&
\dfrac{\partial (u_\lambda^{(1)} \, - \, u_\lambda^{(2)})}{\partial t} \, - \,
 \triangle ( u_\lambda^{(1)} \, - \, u_\lambda^{(2)} )
\\
%3
&
  \, + \, \lambda^{1-\kappa} \bigl[
   \dot \chi \bigl( \bigl( | u_\lambda^{(1)} |^2 \, - \, 1 \bigr)^2 \bigr)
    \bigl( | u_\lambda^{(1)} |^2 \, - \, 1 \bigr) u_\lambda^{(1)} 
     \, - \, 
      \dot \chi \bigl( \bigl( | u_\lambda^{(2)} |^2 \, - \, 1
       \bigr)^2 \bigr)
        \bigl( | u_\lambda^{(2)} |^2 \, - \, 1 \bigr) u_\lambda^{(2)}
         ]
\\
%4
&
\; = \; 0
 \quad \mathrm{in} \quad Q(T),
\\[2mm]
%2
&
u_\lambda^{(1)} \, - \, u_\lambda^{(2)}
 \; = \; u_0^{(1)} \, - \, u_0^{(2)}
  \quad \mathrm{on} \quad \partial Q(T),
\end{array}
\right.
\label{EQ:GLHF-Diff}
\end{equation}
by applying Theorem 9.1 in 
 Lady\v{z}henskaya,~O.~A., Solonnikov,~V.~A., Ural'ceva,~N.~N.
  \cite[p.341]{ladyzhenskaya-solonnikov-uralceva}
   to \eqref{EQ:GLHF-Diff}, 
    we assert \eqref{INEQ:Schauder}.
$\qed$
\begin{Rem}\label{REM:GLHF-Classical}
We call the classical solution in Theorem \ref{THM:Existence-GLHF}
 \lq\lq the classical GLHF.\rq\rq
\end{Rem}
\par
Next we introduce a crude bound:
\begin{Lem}{\rm{(Crude Bound).}}\label{LEM:CB}
A parabolic analogue to 
 Bethuel,~F., Brezis,~H. and H\'{e}lein,~R. {\rm{\cite[Lemma A.1]{bethuel-brezis-helein}}}
  tells us that
\begin{align}
%1
&
|| \nabla u_\lambda ||_{L_{\mathrm{loc}}^\infty (Q(T))}
 \; \le \; \frac C{\sqrt{\lambda}},
  \quad
   || \nabla^2 u_\lambda ||_{L_{\mathrm{loc}}^\infty (Q(T))}
    \; \le \; \frac C{\lambda},
\notag
\\
%2
&
|| \nabla^3 u_\lambda ||_{L_{\mathrm{loc}}^\infty (Q(T))}
 \; \le \; \frac C{\lambda \sqrt{\lambda}},
  \quad
   || \nabla^4 u_\lambda ||_{L_{\mathrm{loc}}^\infty (Q(T))}
    \; \le \; \frac C{\lambda^2}
%     \quad
%      || \nabla^5 u_\lambda ||_{L^\infty}    
%       \; \le \; \frac C{\lambda^2 \, \sqrt{\lambda}}.
\label{INEQ:CB}
\end{align}
hold for the classical GLHF.
\end{Lem}
Next we prove a maximal principle:
\begin{Thm}{\rm{(Maximal Principle).}}\label{THM:Max}
Each of the classical GLHF $\{ u_\lambda \}$ $( \lambda > 0 )$ satisfies
\begin{equation}
| u_\lambda | \; \le \; 1
 \quad \mathrm{for \enspace any \enspace point}
  \; z \, \in \, Q(T).
\label{INEQ:Max}
\end{equation}
\end{Thm}
\vskip 9pt
\noindent{\underbar{{Proof of Theorem \ref{THM:Max}}}.}
\rm\enspace
\vskip 6pt
Set the truncation function $( |u_\lambda|^2 -1 )^{(0)}$ as
\allowdisplaybreaks\begin{align*}
%1
&
( |u_\lambda|^2 -1 )^{(0)} \; = \; 
\begin{cases}
%1.1
0 & (|u_\lambda|^2 \; \le \; 1),
\\
%1.2
|u_\lambda|^2 -1 & (|u_\lambda|^2 \; > \; 1).
\end{cases}
\end{align*}
A multiplier of \eqref{EQ:GLHF} by 
 $(|u_\lambda|^2 -1)^{(0)} u_\lambda$ 
  and integrate it on $(0,t)$ $\times$ $\mathbb{B}^d$
   with any $t$ in $(0,T)$ observes
\allowdisplaybreaks\begin{align}
%1
&
\lint_0^t \, dt \lint_{\mathbb{B}^d}
 \Bigl\la \frac{\partial u_\lambda}{\partial t}, 
  u_\lambda \Bigr\ra 
   (|u_\lambda|^2 -1 )^{(0)} \, dx
\label{EQ:Max-1}\\
%2
&
\, + \, \lint_0^t \, dt \lint_{\mathbb{B}^d}
 | \nabla u_\lambda |^2 (|u_\lambda|^2 -1 )^{(0)} \, dx
\, + \, 
 \frac 12 \lint_0^t \, dt \lint_{\mathbb{B}^d} | \nabla (|u_\lambda|^2 -1)^{(0)} |^2 \, dx
\notag\\
%3
&
\, + \, \lint_0^t \lambda^{1-\kappa} \, dt
 \lint_{\{x \in \mathbb{B}^d ; |u_\lambda| \ge 1\}} 
  \dot\chi \bigl( (| u_\lambda |^2 \, - \, 1)^2 \bigr) (| u_\lambda |^2 \, - \, 1) 
   (|u_\lambda|^2 -1)^{(0)} | u_\lambda |^2 \, dx \; = \; 0.
\notag
\end{align}
Note that 
\allowdisplaybreaks\begin{align}
%1
&
\lint_0^t \, dt \lint_{\mathbb{B}^d}
 \left\la \frac{\partial u_\lambda}{\partial t},  u_\lambda \right\ra
  (|u_\lambda|^2 - 1)^{(0)} \, dx
\; = \; \frac 14 \lint_0^t \, dt \frac d {dt} \lint_{\mathbb{B}^d}
 \bigl( ( | u_\lambda |^2 \, - \, 1 )^{(0)} \bigr)^2 \, dx.
\notag
%\label{INEQ:Max-1}
\end{align}
Since the second, the third and the fourth terms in \eqref{EQ:Max-1} are
 nonnegative, we thus infer
\begin{gather}
\lint_0^t dt \,
 \frac d {dt} \lint_{\mathbb{B}^d}
  \bigl( ( | u_\lambda |^2 \, - \, 1 )^{(0)} \bigr)^2 \, dx
   \; \le \; 0.
\label{INEQ:Max-2}
\end{gather}
\par
So we arrive at
\begin{equation*}
\lint_{\mathbb{B}^d}
 \bigl( ( | u_\lambda (t) |^2 \, - \, 1 )^{(0)} \bigr)^2 \, dx
  \; \le \;
   \lint_{\mathbb{B}^d}
    \bigl( ( | u_0 |^2 \, - \, 1 )^{(0)} \bigr)^2 \, dx,
\end{equation*}
which can read $|u_\lambda| \le 1$ at almost all $z$ $\in$ $Q(T)$.
 Since the mapping $u_\lambda$ is continuous, we claim \eqref{INEQ:Max}.
  \qed
%%
%\begin{Rem}{}
%\par\noindent
%Since $\dot\chi\bigl((|u_\lambda|^2 - 1)^2\bigr)$
% $=$ $1$ 
%  in $0 \le |u_\lambda| \le 1$,
%   we can write \eqref{EQ:GLHF} as
%%
%\begin{equation}
%\frac{\partial u_\lambda}{\partial t}
% \, - \, \triangle u_\lambda
%  \, + \, \lambda^{1-\kappa}
%   \bigl( | u_\lambda |^2 \, - \, 1 \bigr) u_\lambda \; = \; 0.
%\label{EQ:GL}
%\end{equation}
%\end{Rem}
%
\par
We mention three fundamental energy inequalities.
 We only mimic those for the usual linear heat flow:
Then a multiplier of \eqref{EQ:GLHF} 
 by ${\partial u_\lambda}/{\partial t}$,
  an integration of it on $Q(T)$
   permit us to state
%
%%%%%%%%%%%%%%%%%%%%%%%%%%%%%%%%%%%%%%%%%%%%%%%%%%%%%%%%%%%%%%%%%%%%%%%
%
\par
\begin{Thm}{\rm{(Energy Estimate).}}\label{THM:Energy-Estimate}
For any numbers $t_1$ and $t_2$ with 
 $0$ $\le$ $t_1$ $\le$ $t_2$ $\le$ $T$, 
  the classical GLHF $u_\lambda$ satisfies
\allowdisplaybreaks\begin{align}
%1
\lint_{t_1}^{t_2} & \, dt \, \lint_{\mathbb{B}^d}
 \left| \frac{\partial u_\lambda}{\partial t} \right|^2 \, dx
  \, + \, \lint_{\mathbb{B}^d} \mathbf{e}_\lambda (t_2,x) \, dx
\label{INEQ:Energy-Estimate1}
\\
%2
&
\, + \, \frac {\log \lambda}{4}
 \lint_{t_1}^{t_2} \dot{\kappa} \, \lambda^{1-\kappa} \, dt
  \lint_{\mathbb{B}^d}
   \chi (\bigl( | u_\lambda |^2 \, - \, 1 \bigr)^2) \, dx
\; =  \; \lint_{\mathbb{B}^d} \mathbf{e}_\lambda (t_1,x) \, dx.
\notag
\end{align}
\end{Thm}
The third theorem is used for the proof of Theorem \ref{THM:Main-2}:
\begin{Thm}{\rm{(Energy Decay Estimate).}}\label{THM:EDE}
Assume $\left. u_0 \right\vert_{\partial \mathbb{B}^d}$ $=$
 $\mathit{a \; constant}$.
  Then the following 
\allowdisplaybreaks\begin{align}
%1
&
\lint_{\mathbb{B}^d} \mathbf{e}_\lambda (|x|^2+1) \, dx \, e^{(d-2)t}
 \, \le \, \frac 12 \lint_{\mathbb{B}^d} | \nabla u_0 |^2 (|x|^2+1) \,
 dx
\label{INEQ:EDE}
\end{align}
is valid for the classical GLHF and in any time $t$ $\in$ $(0,T)$.
\end{Thm}
\vskip 9pt
\noindent{\underbar{{Proof of Theorem \ref{THM:EDE}}}}.
\vskip 6pt
Multiply \eqref{EQ:GLHF} by $( \partial u_\lambda / \partial t) ( |x|^2 +1 )$
 and integrate it over $\mathbb{B}^d$ to verify
\allowdisplaybreaks\begin{align}
%1
& \lint_{\mathbb{B}^d}
 \left| \frac{\partial u_\lambda}{\partial t} \right|^2
  (|x|^2+1) \, dx
\, + \, \frac d{dt} \lint_{\mathbb{B}^d}
 \mathbf{e}_\lambda (|x|^2+1) \, dx
\label{EQ:EDE-1}
\\
%2
&
\, + \, 2 \lint_{\mathbb{B}^d}
 \Bigl\langle 
  x\!\cdot\!\nabla u_\lambda,
   \frac{\partial u_\lambda}{\partial t}
    \Bigr\rangle \, dx
\, + \, {\log \lambda} 
 \frac {\dot{\kappa} \, \lambda^{1-\kappa}}{4} \lint_{\mathbb{B}^d} 
  \chi \bigl( (| u_\lambda |^2 \, - \, 1 )^2 \bigr) (|x|^2+1) \, dx
   \; = \; 0.
\notag
\end{align}
\par
Next, a multiplier of \eqref{EQ:GLHF} by $-2x\!\cdot\!\nabla u_\lambda$ and
 an integration of it over $\mathbb{B}^d$ imply 
\allowdisplaybreaks\begin{align}
%1
& -2 \lint_{\mathbb{B}^d}
 \Bigl\langle \frac{\partial u_\lambda}{\partial t}, x \cdot \nabla u_\lambda \Bigr\rangle \, dx
\, + \, 2(d-2)
 \lint_{\mathbb{B}^d} \mathbf{e}_\lambda \, dx
\label{EQ:EDE-2}
\\
%2
&
\, + \, 
 \lint_{\partial \mathbb{B}^d}
  \biggl( 
   \left| \frac {\partial u_\lambda}{\partial |x|} \right|^2 
    \, - \, 
     | \nabla_{\tau} u_\lambda |^2
      \biggr) \, d {\cal{H}}_{{x}}^{\mathrm{d-1}}
       \; \le \; 0.
\notag
\end{align}
\par
Summing up \eqref{EQ:EDE-1} and  \eqref{EQ:EDE-2},
 noting $\left. u_0 \right\vert_{\partial \mathbb{B}^d}$ $=$
  $\mathit{a \; constant}$ in the trace sense
   and multiplying it by $ e^{(d-2)t}$,
    we arrive at
\begin{equation}
\frac d{dt} \, \biggl( \lint_{\mathbb{B}^d} 
 \mathbf{e}_\lambda (|x|^2+1) \, dx
  \, e^{(d-2)t} \biggr) 
   \; \le \; 0,
\label{INEQ:Last-1}
\end{equation} 
which concludes our result by integrating from $0$ to 
 any positive number $t$ $\in$ $(0,T)$ with respect to $t$.
  $\qed$
\par
By combining Theorem \ref{THM:Energy-Estimate} 
 with the proof of Theorem \ref{THM:EDE},
  we obtain the following inequality.
   We call it \lq\lq A parabolic Pokhojaev inequality\rq\rq.
\vskip 9pt
\begin{Cor}{\rm{(Parabolic Pokhojaev Inequality).}}\label{COR:PPI}
Let the mapping $u_\lambda$ be the classical GLHF{\rm{;}}
 We infer
\allowdisplaybreaks\begin{align}
%1
&
\lint_0^T \, dt
 \lint_{\partial \mathbb{B}^d}
  \left| \frac {\partial u_\lambda}{\partial |x|} \right|^2
   \, d {\cal{H}}_{{x}}^{\mathrm{d-1}}
\notag\\
%2
&
\; \le \; T 
 \lint_{\partial \mathbb{B}^d}
  | \nabla_{\tau} u_0 |^2
    \, d {\cal{H}}_{{x}}^{\mathrm{d-1}}
\, + \, \frac 12 \lint_{\mathbb{B}^d}
 | \nabla u_0 |^2 
  ( |x|^2 + 1 ) \, dx.
\label{INEQ:PPI}
\end{align}
\end{Cor}
\vskip 6pt \par
We finally introduce a local energy inequality without a proof.
 It will be used a several times in the rest of the paper.
%
%%%%%%%%%%%%%%%%%%%%%%%%%%%%%%%%%%%%%%%%%%%%%%%%%%%%%%%%%%%%%%%%%%%%%%%
%
\par
\begin{Thm}{\rm{(Local Energy Inequality).}}\label{THM:LEI}
The following inequality 
\allowdisplaybreaks\begin{align}
%1
& \lint_{P_R (z_0)}
 \left| \frac{\partial u_\lambda}{\partial t} (z) \right|^2 \, dz
\, + \, \Esssup_{t_0 - R^2 < t < t_0 + R^2 }
  \lint_{B_R (z_0)} \mathbf{e}_\lambda (t,x) \, dx
\label{INEQ:LEI}
\\
%2
& 
\; \le \; \frac {C}{R^2} \lint_{P_{2R} (z_0)}
 \mathbf{e}_\lambda (z) \, dz
\notag
\end{align}
holds for the classical GLHF 
 and any parabolic cylinder $P_{2R} (z_0)$
  compactly contained in $Q(T)$.
\end{Thm}
\par
%We claim that the same result in 
Theorem \ref{THM:Existence-GLHF},
  Lemma \ref{LEM:CB}, Theorem \ref{THM:Max}, Theorem \ref{THM:Energy-Estimate},
   Theorem \ref{THM:EDE}, Corollary \ref{COR:PPI} and Theorem \ref{THM:LEI}
are valid for less stringent smoothness requirement for the initial and boundary condition 
  $u_0$ belonging to $H^{1,2} ( \mathbb{B}^d ; \mathbb{S}^D )$
   $\cap$ $H^{2,2} ( \mathbb{B}^d\setminus {B}_{1-\delta_0} (0) ; \mathbb{S}^D )$
    with a positive number sufficiently small $\delta_0$:
Indeed, take the mollifier of the mapping $u_0$ 
 and passing to the limit, we readily see
\begin{Thm}{\rm{(GLHF)}.}\label{THM:GLHF}
Give a mapping $u_0 \, \in \, H^{1,2} ( \mathbb{B}^d \, ; \, \mathbb{S}^D )$
 $\cap$ $H^{2,2} ( \mathbb{B}^d\setminus {{B}_{1-\delta_0} (0)} ; \mathbb{S}^D )$ 
  with a positive number $\delta_0$ sufficiently small,
the GLHF $u_\lambda$, i.e. the mapping 
 {\rm{(a),(b),(c),(d)}} in p.p \pageref{DEF:GLHF} exists.
In addition, the GLHF satisfies 
 Lemma\ref{LEM:CB}, Theorem \ref{THM:Max}, Theorem \ref{THM:Energy-Estimate}
  Theorem \ref{THM:EDE}, Corollary \ref{COR:PPI} and
   Theorem \ref{THM:LEI}.
\end{Thm}
\begin{Rem}{}\label{REM:EQ}
\par\noindent
Since $\dot\chi\bigl((|u_\lambda|^2 - 1)^2\bigr)$ $=$ $1$ 
 in $0 \le |u_\lambda| \le 1$,
  \eqref{EQ:GLHF} and $\mathbf{e}_\lambda$ reduce to
\allowdisplaybreaks\begin{align}
%1
&  
\frac{\partial u_\lambda}{\partial t}
 \, - \, \triangle u_\lambda
  \, + \, \lambda^{1-\kappa}
   \bigl( | u_\lambda |^2 \, - \, 1 \bigr) u_\lambda \; = \; 0,
\label{EQ:GL}
\\
%2
&
\frac 12 |\nabla u_\lambda|^2 \, + \, 
 \frac{\lambda^{1-\kappa}}4
  \chi ((| u_\lambda |^2 \, - \, 1 )^2).
\label{EQ:GL-Energy}
\end{align}
\end{Rem}
\begin{Rem}{}\label{REM:Inequalities}
\par\noindent
Henceforth, if we quote to the one of inequalities above,
 it means the one for the GLHF 
  subject to the initial and boundary mapping
   $u_0 \, \in \, H^{1,2} ( \mathbb{B}^d ; \mathbb{S}^D )$
    $\cap$ $H^{2,2} ( \mathbb{B}^d\setminus {B}_{1-\delta_0} (0) ; \mathbb{S}^D )$
     with a positive number sufficiently small $\delta_0$.
For instance, if we say \lq\lq from \eqref{INEQ:Max} in Theorem
 \ref{INEQ:Max},\rq\rq
  it indicates
   \lq\lq from \eqref{INEQ:Max} in Theorem \ref{INEQ:Max}
    which holds for the GLHF above,\rq\rq
\end{Rem}
%
% glhf_mono
%
\subsection{\enspace Monotonicity For Scaled Energy.}
We introduce a monotonicity inequality for the scaled energy.
 In our settings, see Y.Chen and F.H.Lin \cite{chen-lin} about it's proof.
Likewise Theorem \ref{THM:Energy-Estimate} etc,
 approximate smoothly the boundary condition,
  note \eqref{INEQ:PPI} of Corollary \ref{COR:PPI} and
   pass to the limit to read
% may be originated in Y.Chen and M.Struwe \cite{chen-struwe}.
%
\begin{Thm}{\rm{(Monotonicity for Scaled Energy).}}\label{THM:Mon}
For any point $z_0 \, = \, (t_0,x_0)$ $\, \in \, Q (T)$ and any positive number $R$
 with $t_0  - (2R)^2 > 0$, the scaled energy is denoted by
\allowdisplaybreaks\begin{align}
%1
&
E_\lambda (R;z_0) \; = \;
 \frac 1 {R^d}
  \lint_{t_0 - (2R)^2}^{t_0 - R^2} \, dt
   \lint_{\mathbb{B}^d} \mathbf{e}_\lambda
\, \exp \Bigl( \frac {|x-x_0|^2}{4(t-t_0)}\Bigr) \, dx.
\label{EQ:Res-Ene-Lar}
\end{align}
Then we have
\allowdisplaybreaks\begin{align}
%1
& \qquad
\frac {d E_\lambda}{dR} (R;z_0) \; \ge \; 
\label{INEQ:Mon}
\\
%2
&
\, - \, 
 \frac 1 {R^{d-1}}
  \lint_{t_0 - (2R)^2}^{t_0 - R^2} \frac {t \, - \, t_0}
   {R^2} \, dt \lint_{\mathbb{B}^d }
    \left| \frac {\partial u_\lambda}{\partial t} 
     \; + \;
      \frac {x \, - \, x_0}{2(t \, - \, t_0)} \cdot \nabla u_\lambda \right|^2
       \exp \Bigl( \frac{|x \, - \, x_0|^2}{4(t \, - \, t_0)} \Bigr) \, dx
\notag\\
%3
&
\; + \; \frac 1 {2R^{d+1}}
 \lint_{t_0 - (2R)^2}^{t_0 - R^2} 
  \lambda^{1-\kappa}
   \lint_{\mathbb{B}^d }
    \bigl( | u_\lambda |^2 \, - \, 1 \bigr)^2
     \exp \Bigl( \frac{|x \, - \, x_0|}{4(t \, - \, t_0)} \Bigr)
      \, dx
\; - \;  C_{\mathsf{M}} R
\notag
\\
%6
&
\text{with} \; d_0 \; = \; {d} \, ( x_0, \partial \mathbb{B}^d )
 \; \mathit{and} \; 
  C_{\mathsf{M}} \; = \; \frac C {d_0^{d+2}}
   \biggl(
    \lint_{\partial \mathbb{B}^d}| \nabla_\tau u_0 |^2 \, d\mathcal{H}_x^{d-1}
     \, + \,
      \lint_{\mathbb{B}^d} | \nabla u_0 |^2 \, dx \biggr).
\notag
\end{align}
\end{Thm}
\begin{Cor}\label{COR:Mon}
An integration of \eqref{INEQ:Mon} from $R_1$ to $R_2$ with 
 $0 < R_1 < R_2$ over $R$ yields
\allowdisplaybreaks\begin{align}
%1
E_\lambda & (R_2;z_0) \; \ge \; E_\lambda (R_1;z_0)
\notag\\
%2
&
\, + \, \lint_{R_1}^{R_2} \frac {dR}{R^{d-1}} \,
 \lint_{t_0 - (2R)^2}^{t_0 - R^2} \, dt \lint_{\mathbb{B}^d }
  \left| \frac {\partial u_\lambda}{\partial t}
  \; + \;
   \frac {x \, - \, x_0}{2(t \, - \, t_0)} \cdot \nabla u_\lambda
		    \right|^2
    \exp \Bigl( \frac{|x \, - \, x_0|^2}{4(t \, - \, t_0)} \Bigr) \, dx
\notag\\
%3
&
\, - \, \frac {C_\mathsf{M}}2 ( R_2^2 \, - \, R_1^2).
\label{INEQ:MonMon}
\end{align}
\end{Cor}
%
% glhf_rp
%
\subsection{\enspace Hybrid type Inequality for GLHF.}
We are in the position to prove an inequality of the hybrid type;
 This inequality is the one of the crucial tools in the paper.
  We claim
\begin{Thm}{\rm{(Hybrid Inequality).}}\label{THM:HI}
For any positive number $\epsilon_0$,
 any point $z_0$ in $Q(T)$, 
  setting $d_0$ as $d_0$ $=$ $\dist (z_0,\partial Q(T))$,
   there exists a positive constant
    $C(\epsilon_0, d_0)$ satisfying $C(\epsilon_0, d_0)$
     $\nearrow$ $\infty$ as $\epsilon_0 \searrow 0$
      or $d_0 \searrow 0$ such that
\allowdisplaybreaks\begin{align}
%1
\lint_{P_{R} (z_0)} \mathbf{e}_\lambda (z) \, dz
 & \, \le \; \epsilon_0 \lint_{P_{2R} (z_0)} \mathbf{e}_\lambda (z) \, dz
  \, + \, \frac {C(\epsilon_0, d_0)}{R^2} \lint_{P_{2R} (z_0)}
   | u_\lambda (z) \, - \, a (t) |^2 \, dz
\notag
\\
%2
& \, + \, \frac {C (\epsilon_0,R, d_0)}{\log \lambda}, 
\label{INEQ:HI}
\end{align}
where $P_{2R} (z_0)$ is any parabolic cylinder compactly contained
 in $Q(T)$ and $a$ $=$ $a(t)$ $=$ $(a^i (t))$
  $(i = 1,2, \ldots, D+1)$
   is any $L^2$-mapping with respect to a positive parameter $t$.
\end{Thm}
\par
As a preliminary we list symbols and auxiliary functions employed only here.
 Give $L_\lambda$ by $[ \lambda^{3/(1+\theta_0)} ]$.
  We then  first introduce the decomposition convention:
   Put
\allowdisplaybreaks\begin{align}
%1
&
\triangle \theta_l \; = \;  (1/2)^l, \; 
 \triangle r_l \; = \; \triangle \rho_l
  \; = \; r (1/2)^l
   \quad ( l = 1,2, \ldots, L_\lambda ),
\notag
\\
%3
\rho_l \; = \; & 
\left\{
\begin{array}{ll}
0
 & \quad ( l = 0 )
\\[2pt]
( 1 - \epsilon_0^4 ) r
 & \quad ( l = 1 )
\\[2pt]
( 1 - \epsilon_0^4 ) r
 \, + \, 
  C_1 \epsilon_0^4 r
   \sum_{j=1}^{l-1}
    \triangle \theta_j
& \quad ( l = 2, \ldots, L_\lambda )
\end{array}
\right.
\notag
\\
%5
&
\textrm{with} \;
 C_1 \, = \, 
  \bigl(\sum_{l=1}^{L_\lambda-1} \triangle \theta_l \bigr)^{-1},
\notag
\\
%6
&
\left\{ 
\begin{array}{l}
 \rho_l^{\frak{b}} \; = \; 
  \rho_l \, - \, 
   2 \triangle \rho_l/3
\\
\rho_l^{\frak{f}} \; = \;
 \rho_l \, - \,
  \triangle \rho_l/3
\end{array}
\right.
(l \, = \, 1,2,\ldots,L_\lambda).
\notag
\end{align}
Throughout $t$ $\in$ $(-r^2,r^2)$,
 choose numbers $r_l^{\frak{b}}$, $r_l^{\frak{f}}$ $( l \, = \, 1,2,\ldots, L_\lambda )$ 
  so that they satisfy
\begin{align}
%1
&
\frac{\triangle \rho_l}{12}
 \lint_{\{\rho_l^{\frak{b}}\} \times \mathbb{S}^{d-1}}
  | \nabla u_\lambda ( t,x ) |^2 \, %
    d\mathcal{H}_x^{d-1}
\; = \; %\kern-15mm
 \lint_{\rho_{l}^{\frak{b}} - \triangle \rho_l/12}
      ^{\rho_{l}^{\frak{b}} + \triangle \rho_l/12}
   \rho^{d-1} \, d\rho
    \lint_{\mathbb{S}^{d-1}}
     | \nabla u_\lambda ( t,x ) |^2 \, d\omega_{d-1},
\notag
\\
%2
&
\frac{\triangle \rho_l}{12}
 \lint_{\{\rho_l^{\frak{f}}\} \times \mathbb{S}^{d-1}}
  | \nabla u_\lambda ( t,x ) |^2 \, d\mathcal{H}_x^{d-1}
\; = \; 
 \lint_{\rho_{l}^{\frak{f}} - \triangle \rho_l/{12}}
  ^{\rho_{l}^{\frak{f}} + \triangle \rho_l/{12}}
   \rho^{d-1} \, d\rho
    \lint_{\mathbb{S}^{d-1}}
     | \nabla u_\lambda ( t,x ) |^2 \, d\omega_{d-1}
\notag
\\
%3
&
\mathrm{and} \quad 
 r_1^{\frak{f}} \, = \, \rho_1,
  r_{0} \, = \, 0, \;
   r_{L_\lambda+1}^{\frak{b}} \, = \, r,
\notag
\\
%3
\widetilde{r}_l \; = \; &
\left\{
\begin{array}{ll}
0
 & \quad ( l = 0 )
\\[2pt]
( 1 - \epsilon_0^4 ) r
 & \quad ( l = 1 )
\\[2pt]
( 1 - \epsilon_0^4 ) r
 \, + \,
  C_2 \epsilon_0^4 r
   \sum_{j=1}^{l-1}
    \triangle \theta_j^3
& \quad ( l = 2, \ldots, L_\lambda )
\notag
\end{array}
\right.
\notag
\\
%4
&
\textrm{with} \;
 C_2 \, = \,
  \bigl(\sum_{l=1}^{L_\lambda-1 } \triangle \theta_l^3 \bigr)^{-2}.
\notag
\end{align}
\par
Next introduce a mapping $f_\lambda$ which is the solution to
\allowdisplaybreaks\begin{align}
%1
&
\left\{
\begin{array}{rcl}
\nabla_{|x|} f_\lambda                                                                                                                                                            
 \, + \, \dfrac {r}{2(d-1)}                                                                                                                                                       
  \triangle_\tau f_\lambda \; & = \; 0                                                                                                                                            
   & \quad \mathrm{in} \quad [ 0, {r}) \times                                                                                                                                     
    \mathbb{S}^{d-1}
\\
%2                                                                                                                                                                                
f_\lambda \; & = \; u_\lambda                                                                                                                                                     
 & \quad \mathrm{on} \quad \{ {r} \} \times                                                                                                                                       
  \mathbb{S}^{d-1}.                                                                                                                                                               
\end{array}                                                                                                                                                                       
\right.                                                                                                                                                                           
\label{EQ:Support}                                                                                                                                                            
\end{align}
\par%
Designate four sorts of annulus
\allowdisplaybreaks\begin{align}
%1
&
T_l^\prime
 \; = \;
  [ \, {\rho}_{l-1}, {\rho}_{l}) \, \times \, {\mathbb{S}^{d-1}},
\quad
T_l^{\frak{b}}
 \; = \;
  [ {r}_{l}^{\frak{b}}, {r}_{l}^{\frak{f}}) \, \times \, {\mathbb{S}^{d-1}},
\notag
\\
%2
&
T_l^{\frak{f}}
 \; = \;
  [ {r}_{l}^{\frak{f}}, {r}_{l+1}^{\frak{b}}) \, \times \, {\mathbb{S}^{d-1}},
\quad
T_l \; = \; T_l^{\frak{b}} \cup T_l^{\frak{f}}
 \; = \; 
  [ \, {r}_{l}^{\frak{b}}, {r}_{l+1}^{\frak{b}}) \, \times \, {\mathbb{S}^{d-1}}
\notag
\\
%3
&
( l = 1,2, \ldots, L_\lambda ).
\end{align}
\par
A first step to prove a Hybrid type inequality is to
 inductively construct a certain support mappings
  ${w}_{\lambda,l}^{\prime}$, 
   ${w}_{\lambda,l}^{\frak{b}}$ and ${w}_{\lambda,l}^{\frak{f}}$
    by making the best of our support function $f_\lambda$:
They are the solutions of
\allowdisplaybreaks\begin{align}
%1
&
\left\{
\begin{array}{rcl}
%1.1
&
- \triangle {w}_{\lambda,l}^{\prime} \; = \; 0 
& 
\quad \mathrm{in} \quad T_l
\\
%1.2
&
\left.
 {w}_{\lambda,l}^{\prime} 
  \right\vert_{\{ \rho_{l-1}\} \times {\mathbb{S}^{d-1}}}
   \; = \;
&
\left. f_\lambda \right\vert_{\{ \tilde{r}_{l-1}\} \times {\mathbb{S}^{d-1}}}
\\
%1.3
&
\left.
 {w}_{\lambda,l}^{\prime}
  \right\vert_{\{ \rho_{l}\} \times {\mathbb{S}^{d-1}}}
   \; = \;
&
\left. f_\lambda \right\vert_{\{ \tilde{r}_{l}\} \times {\mathbb{S}^{d-1}}}
\end{array}
\right.
( l = 2, \ldots, L_\lambda ),
\label{EQ:Support-1}
\\
%2
&
\left\{
\begin{array}{rcl}
%2.1
&
- \triangle {w}_{\lambda,1}^{\prime} \; = \; 0 
& 
\quad \mathrm{in} \quad [0,r_1^{\frak{f}}) \times {\mathbb{S}^{d-1}}
\\
%2.2
&
\left. {w}_{\lambda,l}^{\prime} 
 \right\vert_{\{r_1^{\frak{f}}\} \times {\mathbb{S}^{d-1}}}
  \; = \; \left. f_\lambda
   \right\vert_{\{\tilde{r}_{1}\} \times {\mathbb{S}^{d-1}}},
&
%\quad \mathrm{on} \quad \{r_1^{\frak{f}}\} \times {\mathbb{S}^{d-1}},
\end{array}
\right.
\label{EQ:Support-2}
\\
%3
&
\left\{
\begin{array}{rcl}
%3.1
&
- \triangle {w}_{\lambda,l}^{\frak{b}} \; = \; 0 
&
\quad \mathrm{in} \quad T_l^{\frak{b}}
\\
%3.2
&
\left.
 {w}_{\lambda,l}^{\frak{b}} 
  \right\vert_{\{ r_{l}^{\frak{b}}\} \times {\mathbb{S}^{d-1}}}
   \; = \;
&
 \left. {w}_{\lambda,l}^{\prime}
  \right\vert_{\{ r_{l}^{\frak{b}}\} \times {\mathbb{S}^{d-1}}}
\\
%3.3
&
\left.
 {w}_{\lambda,l}^{\frak{b}}
  \right\vert_{\{ r_{l}^{\frak{f}}\} \times {\mathbb{S}^{d-1}}}
   \; = \;
&
 \left. {w}_{\lambda,l}^{\prime}
  \right\vert_{\{ r_{l}^{\frak{f}}\} \times {\mathbb{S}^{d-1}}}
\end{array}
\right.
( l = 2, \ldots, L_\lambda ),
\label{EQ:Support-3}
\\
%4
&
\left\{
{w}_{\lambda,1}^{\frak{b}} \; = \; 
 {w}_{\lambda,1}^{\prime},
\right.
\label{EQ:Support-4}
\\
%5
&
\left\{
\begin{array}{rcl}
%5.1
&
- \triangle {w}_{\lambda,l}^{\frak{f}} \; = \; 0 
&
\quad \mathrm{in} \quad T_l^{\frak{f}}
\\
%5.2
&
\left.
 {w}_{\lambda,l}^{\frak{f}} 
  \right\vert_{\{ r_{l}^{\frak{f}}\} \times {\mathbb{S}^{d-1}}}
   \; = \;
&
 \left. {w}_{\lambda,l}^{\prime}
  \right\vert_{\{ r_{l}^{\frak{f}}\} \times {\mathbb{S}^{d-1}}}
\\
%5.3
&
\left.
 {w}_{\lambda,l}^{\frak{b}}
  \right\vert_{\{ r_{l+1}^{\frak{b}}\} \times {\mathbb{S}^{d-1}}}
   \; = \;
&
 \left. {w}_{\lambda,l}^{\prime}
  \right\vert_{\{ r_{l+1}^{\frak{b}}\} \times {\mathbb{S}^{d-1}}}
\end{array}
\right.
( l = 2, \ldots, L_\lambda )
\label{EQ:Support-5}
\end{align}
and then set mapping ${w}_{\lambda,l}$ by
\allowdisplaybreaks\begin{align}
%1
&
{w}_{\lambda,l} \; = \;
 \left\{
\begin{array}{rcl}
%%1.1
& 
{w}_{\lambda,l}^{\frak{b}}
&
\mathrm{in} \quad
 [ r_l^{\frak{b}}, r_l^{\frak{f}}) \times {\mathbb{S}^{d-1}}
\\
%1.2
&
{w}_{\lambda,l}^{\frak{f}}
&
\mathrm{in} \quad
 [ r_l^{\frak{f}}, r_{l+1}^{\frak{b}}) \times {\mathbb{S}^{d-1}}
\end{array}
\right.
( l = 2, \ldots, L_\lambda ),
\\
%2
&
{w}_{\lambda,1} \; = \; {w}_{\lambda,1}^{\frak{f}}
  \quad \mathrm{in} \quad [0,r_1^{\frak{f}}) \times {\mathbb{S}^{d-1}}.
\end{align}
\par
We state a property on the mappings $w_{\lambda,l}$ used below.
\begin{Lem}\label{LEM:W-1}
The mappings ${w}_{\lambda,l}^{\frak{b}}$ and ${w}_{\lambda,l}^{\frak{f}}$ have the following property{\rm{:}}
\allowdisplaybreaks\begin{align}
%1
&
| x \cdot \nabla w_{\lambda,l}^{\frak{b}} |^2
 \, + \, 
  | x \cdot \nabla w_{\lambda,l}^{\frak{f}} |^2
\notag
\\
%2
&
\; \le \; \frac C{\triangle r_l^{d-2}}
\Bigl(
 \lint_{\widetilde{r}_{l-1}}^{\widetilde{r}_l} \rho^{d-1} d\rho
  \lint_{\{ \rho \} \times {\mathbb{S}^{d-1}} \cap B_{\triangle r/2} (\rho,x/\rho)}
   | \nabla_{\rho} f_\lambda |^2
    \, d\mathcal{H}_y^{d-1}
\notag
\\
%3
\quad &
\, + \, \triangle r_l
 \lint_{\{ \widetilde{r}_l \} \times {\mathbb{S}^{d-1}} \cap B_{\triangle r/2} (\widetilde{r}_l,x/\widetilde{r}_l)}
  | \nabla_{\tau} f_\lambda |^2 \, d\mathcal{H}_y^{d-1}
   \Bigr)
\label{INEQ:W-1}
\end{align}
holds for any point $x$ in $T_l$.

\end{Lem}
\noindent{\underbar{Proof of Lemma \ref{LEM:W-1}.}}
\rm\enspace 
\par\vskip 6pt
The estimates for $ x \cdot \nabla w_{\lambda,l}^{\frak{b}}$ 
 is performed by the straight-forward computation from the explicit
  formula, the mean value theorem for Laplace equation and the sub-harmonic estimates.
\allowdisplaybreaks\begin{align}
%1
&  
w_{\lambda,l}^{\frak{b}} (x)
\; = \; w_{\lambda,l}^{\frak{b}} (\rho, \omega_{d-1}/\rho)
\notag
\\
%2
&
\; = \; \sum_{\substack{n=1,\\ \alpha \in N(n)}}^\infty
 a_n^{(\alpha)} \Bigl( \frac \rho {r_{l}} \Bigr)^n
  \phi_n^{(\alpha)} (\omega_{d-1})
\, + \, 
 \sum_{\substack{n=1,\\ \alpha \in N(n)}}^\infty b_n^{(\alpha)}
  \Bigl( \frac {r_{l-1}} \rho \Bigr)^n
   \phi_n^{(\alpha)} (\omega_{d-1})
\, + \, a_0 \phi_0 (\omega_{d-1})
\notag
\\
\intertext{with}
%3
&
a_n^{(\alpha)} \; = \;
 \frac {f_\lambda^{n, (\alpha)} (t,\tilde{r}_{l-1}) \tau_l^{n}
        \, - \, 
        f_\lambda^{n, (\alpha)} (t,\tilde{r}_{l})}
        {\tau_l^{2n} - 1},
\notag
\\
%4
&
b_n^{(\alpha)} \; = \;
 \frac {f_\lambda^{n, (\alpha)} (t,\tilde{r}_l) \tau_l^{n}
        \, - \,
        f_\lambda^{n, (\alpha)} (t,\tilde{r}_{l-1})}
        {\tau_l^{2n} - 1},
\notag
\\
%4
&
\tau_l \; = \; \frac {r_{l-1}}{r_{l}},
\notag
\\
%5
&
f_\lambda^{n, (\alpha)} (t,r)
 \; = \; 
  \frac 1 {|{\mathbb{S}^{d-1}}|}
   \lint_{{\mathbb{S}^{d-1}}}
    \langle f_\lambda (t,r,\omega_{d-1}),
     \phi_n^{(\alpha)} (\omega_{d-1}) \rangle
      \, d\omega_{d-1},
\end{align}
where $\{ \phi_n^{(\alpha)}\}$ $(n=0,1,\ldots \, ; \, \alpha \in N(n))$ is a sequence of the
 independent hyper-spherical harmonics
  and $N(n)$ is the number of independent hyper-spherical harmonics 
   with degree $n$.
The estimate on $ x \cdot \nabla w_{\lambda,l}^{\frak{f}}$ is similar.
$\qed$
%\input estimate_hhf.tex

%% file: estimate_hhf_arch.tex
 %
%#! platex hhf_arch
%
\vskip 2pt
\par
After the preparation above, 
 we show the proof of Theorem \ref{THM:HI}.
\vskip 9pt
\noindent{\underbar{Proof of Theorem \ref{THM:HI}.}}
\vskip 9pt
\rm\enspace
Take the difference between \eqref{EQ:GL} and 
 $- \triangle w_{\lambda,l}$ $=$ $0$ on
  $T_l^{\frak{b}}$ and $T_l^{\frak{f}}$,
   multiplying it by 
    $-2x \cdot \nabla ( u_\lambda \, - \, w_{\lambda,l})$,
     integrate it on $T_l^{\frak{b}}$ and $T_l^{\frak{f}}$ and sum up it for $l$ to verify
\allowdisplaybreaks\begin{align}
%1
&
-2 \, \sum_{l=1}^{L_\lambda} 
 \lint_{T_l}
\left\langle
 \frac {\partial u_\lambda}{\partial t},
  x \!\cdot\! \nabla ( u_\lambda \, - \, w_{\lambda,l} )
 \right\rangle \,  \, dx
\notag
\\
%2
&
\, + \, (d-2) \, 
 \sum_{l=1}^{L_\lambda} 
   \lint_{T_l}
   | \nabla ( u_\lambda \, - \, w_{\lambda,l} )|^2 \,  \, dx
\, + \, \frac {d \, \lambda^{1-\kappa}}{2} \,
 \sum_{l=1}^{L_\lambda} 
   \lint_{T_l}
   ( | u_\lambda |^2 \, - \, 1 )^2 \,  \, dx
\notag\\
%3
&
\; = \; -2 \lambda^{1-\kappa}
 \sum_{l=1}^{L_\lambda} 
   \lint_{T_l}
   ( | u_\lambda |^2 \, - \, 1 )
    \la u_\lambda, x \cdot \nabla w_{\lambda,l} \ra
     \,  \, dx
\notag\\
%4
&
\, - \,
 \sum_{l=1}^{L_\lambda-1} 
  \,
\Bigl( \, r_{l+1}^{\frak{b}} \,
 \lint_{ \{ r_{l+1}^\frak{b}\} \times {\mathbb{S}^{d-1}}}
  \, - \,
   r_{l}^{\frak{f}} \lint_{ \{ r_{l}^{\frak{f}} \} \times {\mathbb{S}^{d-1}}}
    \Bigr)
\, | \nabla_{|x|} ( u_\lambda \, - \, w_{\lambda,l}^{\frak{f}})|^2
\, d\mathcal{H}_x^{d-1}
\notag\\
%5
&
\, - \,
 \sum_{l=1}^{L_\lambda-1}
  \,
\Bigl( \, r_{l}^{\frak{f}} \,
 \lint_{ \{ r_{l}^\frak{f}\} \times {\mathbb{S}^{d-1}}}
  \, - \,
   r_{l}^{\frak{b}} 
    \lint_{ \{ r_{l}^{\frak{b}} \} \times {\mathbb{S}^{d-1}}}
     \Bigr)
\, | \nabla_{|x|} ( u_\lambda \, - \, w_{\lambda,l}^{\frak{b}})|^2
 \, d\mathcal{H}_x^{d-1}
\notag
\\
%6
&
\, + \,
 \sum_{l=1}^{L_\lambda-1} 
\,
\Bigl(  r_{l+1}^{\frak{b}}
 \lint_{ \{ r_{l+1}^{\frak{b}} \} \times {\mathbb{S}^{d-1}}}
  \, - \, r_{l}^{\frak{b}}
   \lint_{ \{ r_{l}^{\frak{b}}\} \times {\mathbb{S}^{d-1}}}
    \Bigr)
\, 
 | \nabla_\tau ( u_\lambda \, - \, w_{\lambda,l}^{\frak{b}} ) |^2 \, d\mathcal{H}_x^{d-1}
\notag
\\
%7
&
\, + \,  \frac {\lambda^{1-\kappa}}2
 \sum_{l=1}^{L_\lambda-1} 
 \,
\Bigl( r_{l+1}^{\frak{b}}
 \lint_{ \{ r_{l+1}^{\frak{b}} \} \times {\mathbb{S}^{d-1}}}
  \, - \, r_{l}^{\frak{b}}
   \lint_{ \{ r_{l}^{\frak{b}}\} \times {\mathbb{S}^{d-1}}}
    \Bigr)
( | u_\lambda |^2 \, - \, 1 )^2 \, d\mathcal{H}_x^{d-1}
\notag
\\
%9
&
\; = \; \; (\mathrm{\bigroman{1}}) \, + \, (\mathrm{\bigroman{2}}) \, + \,
 \, \cdots \, \, + \, (\mathrm{\bigroman{5}}).
\label{EQ:RP-0}
\end{align}
\par
From now on, we shall estimate the each term of 
 the right-hand side in \eqref{EQ:RP-0}.
First we estimate the first term ({\bigroman{1}}): 
 Choose a sequence of balls $\{ B_{\triangle r_l (x_{\mathbf{i}_l})}\}$
  $(\mathbf{i}_l \, \in \, I_l)$ with
\allowdisplaybreaks\begin{align}
%1
&
T_l \, \subset \, \cup_{\mathbf{i}_l \, \in \, I_l}
 B_{2 \triangle r_l/3}(x_{\mathbf{i}_l}),
\notag
\\
%2
&
\cup_{\mathbf{i}_l \, \in \, I_l} B_{\triangle r_l}(x_{\mathbf{i}_l})
 \, \subset \,
  T_{l-1} \cup T_l \cup T_{l+1}.
\notag
\end{align}
Notice that 
\allowdisplaybreaks\begin{align}
%1
&
\lambda^{1-\kappa}
 \lint_{B_{2\triangle r_l/3}(x_{\mathbf{i}_l})}
   ( 1 \, - \, | u_\lambda |^2 ) \, dx
\; \le \; C 
 \lint_{B_{\triangle r_l}(x_{\mathbf{i}_l})}
  \mathbf{e}_\lambda \, dx
\, + \, \frac C {\triangle r_l} 
 \lint_{B_{\triangle r_l}(x_{\mathbf{i}_l})}
  | \nabla u_\lambda | \, dx
\notag
\\
%2
&
\, + \, 
 \lint_{B_{\triangle r_l}(x_{\mathbf{i}_l})}
  \left\vert 
   \frac {\partial u_\lambda}{\partial t}
    \right\vert \, dx.
%\label{INEQ:L-3-2}
\notag
\end{align}
We then obtain
\allowdisplaybreaks\begin{align}
%1
(\mathrm{\bigroman{1}}) & \; \le \; \lambda^{1-\kappa}
 \sum_{l=2}^{L_\lambda-1} \sum_{\mathbf{i}_l \in I_l}
\lint_{B_{2\triangle r_l/3}(x_{\mathbf{i}_l})}
 ( 1 \, - \, |u_\lambda|^2 ) 
  \, \ | x \cdot \nabla w_{\lambda,l} | \, dx
\notag
\\
%2
&
\, + \, 
 \lambda^{1-\kappa}                                                    
  \sum_{\mathbf{i}_1 \in I_1}
\lint_{B_{2\triangle r_1/3}(x_{\mathbf{i}_1})}
 ( 1 \, - \, |u_\lambda|^2 )
  \, \ | x \cdot \nabla w_{\lambda,1} | \, dx
\notag
\\
%3
&
\, + \, 
 \lambda^{1-\kappa}                                                    
  \sum_{\mathbf{i}_{L_\lambda} \in I_{L_\lambda}} 
\lint_{B_{2\triangle r_{L_\lambda}/3}(x_{\mathbf{i}_{L_\lambda}})}
 ( 1 \, - \, |u_\lambda|^2 )
  \, \ | x \cdot \nabla w_{\lambda,{L_\lambda}} | \, dx
\notag
\\
%4
&
\; \le \; C
 \sum_{l=2}^{L_\lambda-1} \sum_{\mathbf{i}_l \in I_l}
  \lint_{B_{\triangle r_l}(x_{\mathbf{i}_l})}
   \mathbf{e}_\lambda  
    \, \sup_{x \in B_{2\triangle r_l/3} (x_{\mathbf{i}_l})}
    | x \cdot \nabla w_{\lambda,{l}} | \, dx
\notag
\\
%5
&
\, + \,
 \sum_{l=2}^{L_\lambda-1} \sum_{\mathbf{i}_l \in I_l}
  \lint_{B_{\triangle r_l}(x_{\mathbf{i}_l})}
   \left\vert \frac {\partial u_\lambda}{\partial t} \right\vert  
    \sup_{x \in B_{2\triangle r_l/3} (x_{\mathbf{i}_l})}
     | x \cdot \nabla w_{\lambda,l} | \, dx
\notag
\\
%6
&
\, + \, C
 \sum_{l=2}^{L_\lambda-1} \frac 1{r \triangle r_l}
  \sum_{\mathbf{i}_l \in I_l}
   \lint_{B_{\triangle r_l}(x_{\mathbf{i}_l})}
    | \nabla u_\lambda | \,
     \sup_{x \in B_{2\triangle r_l/3} (x_{\mathbf{i}_l})}
      | x \cdot \nabla w_{\lambda,l} | \, dx
\notag
\\
%7
&
\, + \, C
 \sum_{\mathbf{i}_{1} \in I_1}
  \lint_{B_{\triangle r_{1} (x_{\mathbf{i}_1})}}
   \mathbf{e}_\lambda
    \, \sup_{x \in B_{2\triangle r_1/3} (x_{\mathbf{i}_1})}
     | x \cdot \nabla w_{\lambda,{1}} | \, dx
\notag
\\
%8
&
\, + \, C
 \sum_{\mathbf{i}_{1} \in I_1}
  \lint_{B_{\triangle r_{1} (x_{\mathbf{i}_1})}}
   \left\vert \frac {\partial u_\lambda}{\partial t} \right\vert
    \, \sup_{x \in B_{2\triangle r_1/3} (x_{\mathbf{i}_1})}
     | x \cdot \nabla w_{\lambda,{1}} | \, dx
\notag
\\
%9
&
\, + \, C
 \sum_{\mathbf{i}_{1} \in I_1}
  \lint_{B_{\triangle r_{1} (x_{\mathbf{i}_1})}}
   | \nabla u_\lambda |
    \, \sup_{x \in B_{2\triangle r_1/3} (x_{\mathbf{i}_1})}
     | x \cdot \nabla w_{\lambda,{1}} | \, dx
\notag
\\
%10
&
\, + \, C \lambda \triangle r_{L_\lambda} r^d
 \sup_{B_r \setminus B_{L_\lambda}^{\frak{f}}}
  | \nabla_{|x|} w_{\lambda,L_\lambda}^{\frak{f}}|
\notag
\\
%11
&
\; \le \; C 
 \biggl(
  \lint_{B_{r} \setminus B_{(1-\epsilon_0^4)r}}
   \mathbf{e}_\lambda \, dx
    \biggr)^{1/2}
\;
\biggl(
 \sum_{l=2}^{L_\lambda-1} \sum_{\mathbf{i}_l \in I_l}
  \lint_{B_{\triangle r_l}(x_{\mathbf{i}_l})}
   \mathbf{e}_\lambda \, dx
    \, \sup_{x \in B_{2\triangle r_l/3} (x_{\mathbf{i}_l})}
     | x \cdot \nabla w_{\lambda,l} |^2 \,
\biggr)^{1/2}
\notag
\\
%12
&
\, + \, C
\Bigl(
 \lint_{B_{r} \setminus B_{(1-\epsilon_0^4)r}}
  \left\vert \frac {\partial u_\lambda}{\partial t} \right\vert^2
   \, dx \Bigr)^{1/2} \,
\Bigl(
 \sum_{l=2}^{L_\lambda-1} \sum_{\mathbf{i}_l \in I_l}
  \frac {\triangle r_l^d}{r^2}
   \, \sup_{x \in B_{2\triangle r_l/3} (x_{\mathbf{i}_l})}
    | x \cdot \nabla w_{\lambda,l} |^2
     \Bigr)^{1/2}
\notag
\\
%13
&
\, + \, 
\Bigl(
 \lint_{B_{r} \setminus B_{(1-\epsilon_0^4)r}}
  | \nabla u_\lambda |^2 \, dx \Bigr)^{1/2} \,
\Bigl(
 \sum_{l=2}^{L_\lambda-1} \sum_{\mathbf{i}_l \in I_l}
  \triangle r_l^{d-2}
   \, \sup_{x \in B_{2\triangle r_l/3} (x_{\mathbf{i}_l})}
    | x \cdot \nabla w_{\lambda,l} |^2
     \Bigr)^{1/2}
\notag
\\
%14
&
\, + \, 
 \biggl(
  \lint_{B_{r_1}}
   \mathbf{e}_\lambda \, dx
    \biggr)^{1/2}
\;
\biggl(
 \lint_{B_{r_1}}
  \mathbf{e}_\lambda \, dx
   \, \sup_{x \in B_{r_1}}
    | x \cdot \nabla w_{\lambda,1} |^2 \,
\biggr)^{1/2}
\notag
\\
%15
&
\, + \, C
\Bigl( r^2
 \lint_{B_{r_1}}
  \left\vert \frac {\partial u_\lambda}{\partial t} \right\vert^2
   \, dx \Bigr)^{1/2} \,
\Bigl( r^{d-2}
 \sup_{x \in B_{r_1}}
  | x \cdot \nabla w_{\lambda,1} |^2
   \Bigr)^{1/2}
\notag
\\
%16
&
\, + \, C
\Bigl(
 \lint_{B_{r_1}}
  | \nabla u_\lambda |^2 \, dx \Bigr)^{1/2} \,
\Bigl( r^{d-2}
 \sup_{x \in B_{r_1}}
  | x \cdot \nabla w_{\lambda,1} |^2
   \Bigr)^{1/2}
\notag
\\
%17
&
\, + \, C \lambda^2 \triangle r_{L_\lambda} r^d.
\notag
\end{align}
By applying  Corollary \ref{COR:PPI}, Lemma \ref{LEM:W-1}
 and the sub-harmonic estimate for $x \cdot \nabla w_{\lambda,l}$,
  we can proceed to our evaluation.
\allowdisplaybreaks\begin{align}
%1
(\mathrm{\bigroman{1}}) & 
\; \le \; \frac C {\epsilon_0^2} 
 \lint_{B_{r} \setminus B_{(1-\epsilon_0^4)r}}
  \Bigl(
   \mathbf{e}_\lambda
    \, + \,
     r^2  \left\vert \frac {\partial u_\lambda}{\partial t} \right\vert^2
      \Bigr) \, dx
%\notag
%\\
%%2
%&
\, + \, {C \epsilon_0^2}
 \lint_{B_{r} \setminus B_{(1-\epsilon_0^4)r}}
  | \nabla_{|x|} f_\lambda |^2 \, dx
\notag
\\
%2
&
\, + \, {C \epsilon_0^2}
 \lint_{B_{(1-\epsilon_0^4)r}}
  \bigl(
   \mathbf{e}_\lambda  
    \, + \,
     \left\vert \frac {\partial u_\lambda}{\partial t}
      \right\vert^2
       \bigr) \, dx
\notag
\\
%4
&
\, + \, \frac C{\epsilon_0^2}
 \lint_{B_{(1-\epsilon_0^4)r}}
  | \nabla_{|x|} f_\lambda |^2 \, dx
\, + \, C \lambda^2 \triangle r_{L_\lambda} r^d,
\notag
\end{align}
by the choice of $\triangle r_{L_\lambda}$% $=$ $O(1/\lambda^3)$
 which follows to conclude
\allowdisplaybreaks\begin{align}
%1      
(\mathrm{\bigroman{1}}) & \; \le \; C
\lint_{B_r \setminus B_{r (1-\epsilon_0^4)}}
 \Bigl(
  \mathbf{e}_\lambda
   \, + \,
    r^2  \left\vert \frac {\partial u_\lambda}{\partial t} \right\vert^2
     \Bigr) \, dx
\notag
\\
%2
&
\, + \, C 
 \Bigl(
  \frac 1 {\epsilon_0^2}
   \lint_{B_r \setminus B_{r (1-\epsilon_0^4)}}
    | \nabla_{\tau} u_\lambda |^2 \, dx
     \, + \, 
      \epsilon_0^2 r \lint_{\partial B_r}
       | \nabla_{\tau} u_\lambda |^2 \, d\mathcal{H}_x^{d-1}
        \Bigr)
\label{INEQ:L-1}
\\
%3
&
\, + \, \frac {C(\epsilon_0)}r
 \lint_{\partial B_r}
  | u_\lambda \, - \, a |^2 \, d\mathcal{H}_x^{d-1}
\, + \, \frac {C r^d}{\lambda}.
\notag
\end{align}
Successively this implies
\allowdisplaybreaks\begin{align}
%1
(\mathrm{\bigroman{2}}) & \, + \, (\mathrm{\bigroman{3}})
\; \le \; -r 
  \lint_{\partial B_r} | \nabla_{|x|} u_\lambda |^2 
   \, d\mathcal{H}_x^{d-1}
\notag
\\
%2
&
\, + \,
 2 \sum_{l=1}^{L_\lambda-1} r_{l+1}^{\frak{b}}
  \lint_{ \{ r_{l+1}^{\frak{b}} \} \times {\mathbb{S}^{d-1}}} 
   \langle 
    \nabla_{|x|} u_\lambda, \nabla_{|x|} w_{\lambda,l}^{\frak{f}}
     \rangle \, d\mathcal{H}_x^{d-1}
\notag
\\
%3
&
\, - \, 
 2 \sum_{l=1}^{L_\lambda} r_{l}^{\frak{f}}
  \lint_{ \{ r_l^{\frak{f}} \} \times {\mathbb{S}^{d-1}}}
   \langle
    \nabla_{|x|} u_\lambda, \nabla_{|x|} w_{\lambda,l}^{\frak{f}}
     \rangle \, d\mathcal{H}_x^{d-1}
\notag
\\
%4
&
\, + \, 2 \sum_{l=1}^{L_\lambda} r_{l}^{\frak{f}}
 \lint_{ \{ r_l^{\frak{f}} \} \times {\mathbb{S}^{d-1}}}
  \langle
   \nabla_{|x|} u_\lambda, \nabla_{|x|} w_{\lambda,l}^{\frak{b}}
    \rangle \, d\mathcal{H}_x^{d-1}
\notag
\\
%5
&
\, - \, 2
 \sum_{l=1}^{L_\lambda} r_{l}^{\frak{b}}
  \lint_{ \{ r_l^{\frak{f}} \} \times {\mathbb{S}^{d-1}}}
   \langle
    \nabla_{|x|} u_\lambda, \nabla_{|x|} w_{\lambda,l}^{\frak{b}}
     \rangle \, d\mathcal{H}_x^{d-1}
\notag
\\
%6
&
\, + \, 
 \sum_{l=1}^{L_\lambda} r_{l}^{\frak{f}}
  \lint_{ \{ r_l^{\frak{f}} \} \times {\mathbb{S}^{d-1}}}
   | \nabla_{|x|} w_{\lambda,l}^{\frak{f}} |^2 
    \, d\mathcal{H}_x^{d-1}
%\notag
%\\
%%7
%&
\, + \,
 \sum_{l=1}^{L_\lambda} r_{l}^{\frak{b}}
  \lint_{ \{ r_l^{\frak{b}} \} \times {\mathbb{S}^{d-1}}}
   | \nabla_{|x|} w_{\lambda,l}^{\frak{b}} |^2
    \, d\mathcal{H}_x^{d-1}
\notag
\\
%8
&
\; \le \; - r 
 \lint_{\partial B_r}
  | \nabla_{|x|} u_\lambda |^2 \, d\mathcal{H}_x^{d-1}
\notag
\\
%9
&
\, + \, 
 \sum_{l=1}^{L_\lambda-1} \, \triangle r_l
  \lint_{ \{ r_{l+1}^{\frak{b}} \} \times {\mathbb{S}^{d-1}}}
   | \nabla_{|x|} u_\lambda |^2 
    \, d\mathcal{H}_x^{d-1}
%\notag
%\\
%%10
%&
\, + \,
 \sum_{l=1}^{L_\lambda-1} \, \triangle r_l
  \lint_{ (\{ r_{l}^{\frak{b}} \} \cup \{ r_{l}^{\frak{f}} \}) \times {\mathbb{S}^{d-1}}}
   | \nabla_{|x|} u_{\lambda} |^2 \, d\mathcal{H}_x^{d-1}
\notag
\\
%11
&
\, + \, 
 \sum_{l=1}^{L_\lambda-1} \, \frac {Cr^2}{\triangle r_l}
  \lint_{ \{ r_{l+1}^{\frak{b}} \} \times {\mathbb{S}^{d-1}}}
   | \nabla_{|x|} w_{\lambda,l}^{\frak{f}} |^2
    \, d\mathcal{H}_x^{d-1}
\notag
\\
%12
&
\, + \,
 \sum_{l=1}^{L_\lambda-1} \, \frac {Cr^2}{\triangle r_l}
  \lint_{ (\{ r_{l}^{\frak{b}} \} \cup \{ r_{l}^{\frak{f}} \}) \times {\mathbb{S}^{d-1}}}
   | \nabla_{|x|} w_{\lambda,l}^{\frak{b}} |^2
    \, d\mathcal{H}_x^{d-1}.
\notag
\end{align}
\par
Using Lemma \ref{LEM:W-1}, 
 from a definition on $r_l^{\frak{b}}$ and $r_l^{\frak{f}}$, we arrive at
\allowdisplaybreaks\begin{align}
%1
(\mathrm{\bigroman{2}}) & \, + \, (\mathrm{\bigroman{3}})%
 \; \le \; C
  \Bigl( \frac 1 {\epsilon_0^2}
   \lint_{B_r \setminus B_{(1-\epsilon_0^4)r}} | \nabla u_\lambda |^2 \, dx
\, + \, \epsilon_0^2 r
 \lint_{\partial B_r} | \nabla u_\lambda |^2 
  \, d\mathcal{H}_x^{d-1}
   \Bigr)
\notag
\\
%7
&
\, + \, \frac C {\epsilon_0^2 r} \lint_{\partial B_r}
 | u_\lambda \, - \, a |^2 \, d\mathcal{H}_x^{d-1}
\, + \, \frac {C(\epsilon_0, r)}{\sqrt{\lambda}}.
\label{INEQ:L-2}
\end{align}
\par
From construction on $w_{\lambda,l}$, we readily find that $(\mathrm{\bigroman{4}})$ vanishes.
Finally $(\mathrm{\bigroman{5}})$ becomes
\allowdisplaybreaks\begin{align}
%1
&
(\mathrm{\bigroman{5}}) \; = \;
 \frac{r \lambda^{1-\kappa}}2
  \lint_{\partial B_r}
   ( |u_\lambda|^2 \, - \, 1 )^2 \, d\mathcal{H}_x^{d-1}.
\label{INEQ:L-5}
\end{align}
\par
We also have the following estimate for the left-hand side in \eqref{EQ:RP-0}
  called $(L)$:
\allowdisplaybreaks\begin{align}
%1
\mathrm{(L)} & \, \ge \; \frac {d-2}4
 \lint_{B_{(1-\epsilon_0^4)r}} \mathbf{e}_\lambda \,  \, dx
\, - \, C \, r^2 
 \lint_{B_r}
   \left|
    \frac {\partial u_\lambda}{\partial t}
     \right|^2 \,  \, dx
\label{INEQ:L}
\\
%2
&
\, - \, \frac {C}{r}
 \lint_{\partial B_r}
  | u_\lambda \, - \, a |^2 
   \, d\mathcal{H}_x^{d-1}.
\notag
\end{align}
\par
A substitution of \eqref{INEQ:L-1}, \eqref{INEQ:L-2},
 \eqref{INEQ:L-5} and \eqref{INEQ:L} for \eqref{EQ:RP-0}, 
an integration of it with respect to $t$ 
 $\in$ $(-r^2,r^2)$ verifies
\allowdisplaybreaks\begin{align}
%1
\lint_{P_{r/2}} & \mathbf{e}_\lambda \, dz
\; \le \; C \epsilon_0^2
 \lint_{P_r} 
  \mathbf{e}_\lambda \, dz
\notag
\, + \, C r^2
 \lint_{P_{r}} 
  \biggl| \frac {\partial u_\lambda}{\partial t}  \biggr|^2
   \, dz
\notag
\\
%2
&
\, + \, \frac{C (\epsilon_0)} {r}
 \lint_{-r^2}^{r^2} \, dt \, 
  \lint_{\partial B_{r}}
   | u_\lambda \, - \, a |^2 \, d\mathcal{H}_x^{d-1}
\label{INEQ:RP-2}
\\
%3
&
\, + \, 
 \lint_{-r^2}^{r^2} \frac {r \lambda^{1-\kappa}}2 \, dt  
  \, \lint_{\partial B_{r}}
   ( |u_\lambda|^2 \, - \, 1 )^2 \, d\mathcal{H}_x^{d-1}.
\notag
\end{align}
\par
Integrate \eqref{INEQ:RP-2} from $R/2$ to $R$ with respect to $r$
 and divide it by $R$ to obtain
\allowdisplaybreaks\begin{align}
%1
\lint_{P_{R/4}} & \mathbf{e}_\lambda \, dz
\; \le \; C \epsilon_0
 \lint_{P_{R}} \mathbf{e}_\lambda \, dz
\, + \, 
 \frac {C(\epsilon_0)}{R^2} 
  \lint_{P_{R}} | u_\lambda \, - \, a |^2 \, dz
\notag
\\
%2
&
\, + \, C R^2
 \lint_{P_{R}} \,
  \Bigl| \frac {\partial u_\lambda}{\partial t} \Bigr|^2
   \, dz
\, + \, \frac {C(R,\epsilon_0)}{\log \lambda}.
\label{INEQ:RP-3}
\end{align}
\par
To complete the proof,
 we attempt to cover finely any fixed parabolic cylinder $P_R (z_0)$
  by a family of small parabolic cylinder
   whose diameter is $\epsilon_0 R$.
For this purpose, set $N_1$ be $(2([1/\epsilon_0]+1))^d$.
We equivalently divide the parabolic cube $D_R (z_0)$ $=$ 
 $(t_0-R^2, t_0+R^2)$ $\times$ $C_R (x_0)$
  into small parabolic cubes:
There is a finite sequence of small parabolic cubes
 $\{ D_{\epsilon_0 R} (t_{p'},x_{q'}) \}$ 
  $( p \, = \, 1,2,\ldots, 2([{1/\epsilon_0^2}]+1) \,; \, q \, = \, 1,2,\ldots,N_1)$
   with
\allowdisplaybreaks\begin{align}
%1
&
D_{\epsilon_0 R} (t_{p'},x_{q'}) \, \cap D_{\epsilon_0 R} (t_{r'},x_{s'})
 \; = \; \emptyset \;
  \; \text{if} \; p' \, \ne \, r'
   \; \text{or} \; q' \, \ne \, s',
\notag
\\
%2
&
D_R (z_0) \, \subset \, 
 \bigcup_{p=1}^{2[1/\epsilon_0^2]+2} \bigcup_{q=1}^{N_1}
  \; \mathrm{\lq\lq the} \; \mathrm{closure} \; \mathrm{of} \; {D_{\epsilon_0 R}} (t_{p'},x_{q'})
   \text{\rq\rq}.%{\textquoteright 
\notag
\end{align}
\par
If we take a family of parabolic cylinders $\{ P_{{\epsilon_0} R} (t_p,x_q) \}$
 whose centre $(t_p,x_q)$ is located at all the centre 
  and the vertex of each small parabolic cube $D_{\epsilon_0 R} (t_{p'},x_{q'})$ above, 
   it is fulfilled
\begin{equation}
P_{R} (z_0) \, \subset \, D_{R} (z_0)
 \, \subset \, 
  \bigcup_{(t_p,x_q)} 
   \; \mathrm{\lq\lq the} \; \mathrm{closure} \; \mathrm{of} \;
    P_{{\epsilon_0} R} (t_p,x_q).\text{\rq\rq} %\mathrm{\rq\rq}
\label{COV:1}
\end{equation}
We must remark that
 the number of $P_{4 {\epsilon_0} R} (t_p,x_q)$ that
  includes any fixed point in $D_{R} (z_0)$
   is bounded by $d$.
Once we assign each element of $\{ P_{{\epsilon_0} R} (t_p,x_q) \}$ to $P_{R/4}$
 in \eqref{INEQ:RP-3} and sum it up over all such a parabolic cylinder, we obtain
\allowdisplaybreaks\begin{align}
%1
&
\sum_{p,q}
 \lint_{P_{{\epsilon_0} R} (t_p,x_q)} \mathbf{e}_\lambda \, dz
\; \le \; C \epsilon_0^2
 \sum_{p,q}
  \lint_{P_{4 {\epsilon_0} R} (t_p,x_q)} \mathbf{e}_\lambda \, dz
\notag
\\
%2
&
\, + \, \frac {C(d,\epsilon_0)}{({\epsilon_0} R)^2}
 \sum_{p,q}
  \lint_{P_{4 {\epsilon_0} R} (t_p,x_q)}
   | u_\lambda \, - \, a |^2 \,dz
\, + \, C(d,\epsilon_0) ({\epsilon_0} R)^2
 \sum_{p,q} 
  \lint_{P_{4 {\epsilon_0} R} (t_p,x_q)}
   \left| \frac{\partial u_\lambda}{\partial t} \right|^2 \, dz
\notag
\\
%3
&
\, + \, \frac {C(R,\epsilon_0)}{\log \lambda}.
\label{INEQ:6}
\end{align}
Since the inclusion 
\begin{equation*}
\bigcup_{p,q}
 P_{4 {\epsilon_0} R} (t_p,x_q)
  \; \subset \; 
   P_{3R/2} 
\end{equation*}
holds, it implies that 
 in view of \eqref{INEQ:6} and \eqref{INEQ:LEI} in Theorem \ref{THM:LEI}, 
  taking a positive number $\epsilon_0$ so small with 
   $C \epsilon_0$ $<$ $1$, we describe
\allowdisplaybreaks\begin{align}
%1
&
\lint_{P_{R} (z_0)} \mathbf{e}_\lambda \, dz
\; \le \; \epsilon_0
 \lint_{P_{2R} (z_0)} \mathbf{e}_\lambda \, dz
\, + \, \frac {C ( \epsilon_0 )}{R^2}
 \lint_{P_{2R} (z_0)}
  | u_\lambda \, - \, a |^2 \,dz
\notag
\\
%2
&
\qquad
 \, + \, \frac {C(R,\epsilon_0)}{\log \lambda},
\label{INEQ:RP-5}
\end{align}
which thereby completes the proof. $\qed$

%% file: hhf_hhf_arch.tex
%
%#! platex hhf
%
\setcounter{chapternumber}{3}\setcounter{equation}{0}
\renewcommand{\theequation}%
           {\thechapternumber.\arabic{equation}}
\section{\enspace WHHF}
\subsection{Existence}
This chapter studies the existence and a partial regularity on WHHF.
The existence theorem
 is a slight modification of Y.Chen \cite{chen-89};
  See also L.C.Evans \cite[p.48, 5.A.1]{evans-88} and J.Shatah \cite{shatah}.
First of all we mention convergence theorem directly derived
 from Theorem \ref{THM:Energy-Estimate} in GLHF:
\begin{Thm}{\rm{(Convergence).}}\label{THM:Convergence}
There exist a subsequence $\{ u_{\lambda (\nu)} \}$ $( \nu \, = \, 1,2,\ldots )$ of
 $\{ u_\lambda \}$ $( \lambda > 0 )$ in $V(Q(T) ; \mathbb{S}^D)$
  and a mapping $u \, \in \, $ $V(Q(T);\mathbb{S}^D)$
such that the sequence of mappings $\{ u_{\lambda (\nu)} \}$ 
 $( \nu \, = \, 1,2,\ldots )$
  respectively converges weakly and weakly-$*$ to a mapping $u$ in 
   $H^{1,2} (0,T ; L^2 (\mathbb{B}^d ; \mathbb{R}^{D+1}))$ 
    and $L^\infty (0,T ; H^{1,2} (\mathbb{B}^d ; \mathbb{R}^{D+1}))$.
So does it strongly to the mapping $u$ in
 $L^2 (Q(T) ; \mathbb{S}^D)$ and point-wisely to it
  in almost all $z$ $\in$ $Q(T)$ as $\nu \nearrow \infty$.
\end{Thm}
\par
Theorem \ref{THM:Convergence} enables us state the following existence theorem:
\begin{Thm}{\rm{(Existence).}}\label{THM:Existence-HF}
The GLHF converges to 
 a WHHF in $L^2 (Q(T) ; \mathbb{R}^{D+1})$ 
  as $\lambda \nearrow \infty$ \rm{(}\it{modulo a subsequence of $\lambda$}\rm{)}.
\end{Thm}
\begin{Rem}{}
We hereafter fix the subsequences $\{\lambda ( \nu )\}$ 
 $( \nu=1,2,\ldots )$
  of $\{ \lambda \}$ $( \lambda > 0 )$ 
   chosen in Theorem \ref{THM:Convergence}
\end{Rem}
\subsection{\enspace Partial Regularity}
The next section discusses a partial regularity
 on the WHHF obtained through the limit of the GLHF.
\begin{Def}{}\label{REM:Measure}
Let $\{ u_{\lambda (\nu )}\}$ $( \nu \, = \, 1,2,\ldots )$
 be the sequence selected above and set $\mathbf{e}_{\lambda (\nu)}$
  the Ginzburg-Landau energy density
   $| \nabla u_{\lambda (\nu)} |^2/2$ $+$
    $\lambda (\nu)^{1-\kappa }$ $(| u_{\lambda (\nu)}|^2 \, - \, 1)^2/4$.
We then denotes $\overline{\mathcal{M}}$ by
\begin{equation*}
\overline{\mathcal{M}} (P_R (z_0)) \,= \, 
 \limsup_{\lambda (\nu) \nearrow \infty } \frac 1 {R^d}% \kern -2mm
  \lint_{P_R (z_0)} % \kern-1mm
   \mathbf{e}_{\lambda (\nu) } \, dz
\end{equation*}
where $P_R (z_0)$ is an arbitrary parabolic cylinder compactly contained in $Q(T)$.
\end{Def}
\begin{Lem}{\rm{(Measured Hybrid Inequality).}}\label{LEM:Measure-Hybrid}
Assume that the a sequence of GLHF $\{ u_{\lambda (\nu)} \}$ 
 $( \nu = 1,2,\ldots )$, respectively converges 
  weakly and weakly-$*$ in $H^{1,2} (0,T ; L^2 (\mathbb{B}^d ; \mathbb{R}^{D+1}))$  and 

$L^\infty (0,T;H^{1,2} (\mathbb{B}^d ; \mathbb{R}^{D+1}))$
   to a WHHF $u$ $\, \in \,$ $V( Q ; \mathbb{S}^D )$ as $\lambda (\nu) \nearrow \infty$.
Then take the pass to the limit $\lambda (\nu) \nearrow \infty$ in 
 Theorem \ref{THM:HI} to infer the following{\rm{:}}
For any positive number $\epsilon_0$,
 there exists a positive constant $C ( \epsilon_0 )$
  satisfying $C ( \epsilon_0 )$ $\nearrow \infty$ as $\epsilon_0  \searrow 0$
   such that the inequality
\allowdisplaybreaks\begin{align}
%1
\overline{\mathcal{M}} & \bigl(P_R (z_0) \bigr) 
 \; \le \; 
  \epsilon_0 \overline{\mathcal{M}} \bigl(P_{2R} (z_0) \bigr)
\, + \, C ( \epsilon_0 )
 \avint{P_{2R} (z_0)}  | u (z) \, - \, a (t) |^2 \, dz
\label{INEQ:Measure-Hybrid}
\end{align}
holds for any time-variable $L^2$-mapping $a(t)$ $( = ( a^i (t) )$
 $( i = 1,2, \ldots, D+1 )$
  and any parabolic cylinder $P_{2R} (z_0)$ compactly contained in $Q(T)$.
\end{Lem}
In the similar manner as in L.~Simon \cite[Lemma 2, p.31]{simon95},
 we can assert the following reverse Ponicar\'e inequality.
\begin{Cor}{\rm{(Reverse Poincar\'e Inequality).}}\label{COR:RPI}
Let a mapping $a$ $(= a(t))$ be a time variable mapping in $L^2 ((0,T) \, ;\mathbb{R}^{D+1})$.
 Then \eqref{INEQ:Measure-Hybrid} implies the reverse Poincar\'e inequality
\begin{equation}
R^{d+2} \overline{\mathcal{M}} (P_{R} (z_0))
 \; \le \; C %\kern-3mm
  \lint_{P_{2R} (z_0)} | u (z) \, - \, a (t) |^2 \, dz
\label{INEQ:RPI}
\end{equation}
holds  whenever $P_{2R} (z_0)$ is an arbitrary parabolic cylinder 
 compactly contained in $Q(T)$.
\end{Cor}
%
%%%%%%%%%%%%%%%%%%%%%%%%%%%%%%%%%%%%%%%%%%%%%%%
By combining Corollary \ref{COR:RPI} with Theorem \ref{THM:LEI}
 and using Sobolev imbedding theorem and Poincar\'{e} inequality
  for the space variables,
   we can describe the following lemma.
We refer the proof to Theorem 2.1 in M.~Giaquinta and M.~Struwe \cite{giaquinta-struwe}.
\begin{Lem}{}\label{LEM:RH}
There exists a positive number $q_0$ greater than $1$
 such that differentials $\nabla u$ of the WHHF $u$ belongs to 
  $L_{\loc}^{2q_0}$ $(Q(T) ; \mathbb{R}^{d{(D+1)}})$
   with
\vskip 15pt
\begin{equation}
\biggl( 
 \avint{P_{R} (z_0)} | \nabla u (z) |^{2q_0} \, dz
  \biggr)^{1/q_0}
\; \le \; C
 \avint{P_{2R} (z_0)}
  | \nabla u (z) |^2 \, dz.
\label{INEQ:RH}
\end{equation}
\end{Lem}
If we follow the result by Y.~Chen and M.~Struwe \cite[Lemma 2.4]{chen-struwe},
 we can claim
%%%%%%%%%%%%%%%%%%%%%%%%%%%%%%%%%%%%%%%%%%%%%%%%%%%%%%%%%%
\begin{Thm}{}\label{THM:PR}
For any positive number $\epsilon$, set
\begin{align}
%1
&
\mathbf{sing} ( \epsilon ) \, = \, \bigcap_{R > 0}
 \{ z_0 \, \in \, Q(T) \, ; \,
  \, \overline{\mathcal{M}} \! \bigl(P_R (z_0)\bigr) \, \ge \, \epsilon, \,
   P_R (z_0) \, \subset\!\subset \, Q(T)
\},
\label{EQ:Sing}
\\
&
\mathbf{reg} ( \epsilon ) \; = \; Q(T) \setminus \mathbf{sing}.
\label{EQ:Reg}
\end{align}
Then there exist some positive number $\epsilon_0$ and an increasing function $g(t)$ with 
 $g(0) = 0$ and $g(t)$ $=$ $O (t \log (1/t)^{d+1})$
  $( t \searrow 0 )$ %and the positive number $\epsilon_0$
   such that
    if $z_0$ $\in$ $reg ( \epsilon_0 )$, that is
     for some positive number $R_0$ and positive integer $\lambda_0$
      possibly depending on $z_0$,
\begin{equation}
\frac 1{R_0^d} \lint_{P_{g(R_0)} (z_0)}
 \mathbf{e}_\lambda (z) \, dz
  \; < \; \epsilon_0
\end{equation}
implies
\begin{equation}
\sup_{z \in P_{R_0} (z_0)}
 \mathbf{e}_\lambda (z)
  \; \le \; \frac C{R_0^2},
\label{INEQ:PR}
\end{equation}
as long as any $\lambda$ is more than or equal to
 $\lambda_0$.
\end{Thm}
\begin{Def}{}
In the sequel, we respectively mean $\, \mathbf{sing}$ and $\mathbf{reg}$
  by $\mathbf{sing} ( \epsilon_0 )$ and $\mathbf{reg} ( \epsilon_0 )$, respectively.
\end{Def}
\begin{Lem}{}\label{LEM:Cont-Time}
Pick up any point $z_0$ $\in$ $\mathbf{reg}$
 and fix it.
  On the parabolic cylinder $P_{R_0/2} (z_0)$ which is the half size of
   the cylinder in \eqref{INEQ:PR}, 
    the inequality
\begin{equation}
| u_\lambda (t,x) \, - \, u_\lambda (s,x) |
 \; \le \; \frac C{R_0} |t-s|^{1/2}
\label{INEQ:Cont-Time}
\end{equation}
holds for any points $t$ and $s$ in $[t_0 - (R_0/2)^2, t_0 + (R_0/2)^2]$
 and $x$ $\in$ $\overline{B}_{R_0/2} (x_0)$ with $z_0$ $=$ $(t_0,x_0)$.
\end{Lem}
\noindent{\underbar{Proof of Corollary \ref{LEM:Cont-Time}.}}
\rm\enspace
\vskip 6pt
Assume $s < t$; 
 Then combining Theorem \ref{THM:LEI} with Theorem \ref{THM:PR}
  as $R$ $=$ $\sqrt{t-s}$, we infer
\begin{align}
%1
\lint_{B_{\sqrt{t-s}} (x)} &
 | u_\lambda (t,y) \, - \, u_\lambda (s,y) | \, dy
\notag\\
%2
&
\; \le \; |t-s|^{(d+2)/4}
 \left(
  \lint_s^t \, d\tau %\kern-.5mm
   \lint_{B_{\sqrt{t-s}} (x)}
    \left| \frac {\partial u_\lambda} {\partial \tau} \right|^2
     \, dy \right)^{1/2}
\notag\\
%3
&
\; \le \; C |t-s|^{(d+2)/4}
 \left(
  \frac 1 {|t-s|}
   \lint_{(s+t)/2-(t-s)}^{(s+t)/2+(t-s)} \, d\tau% \kern-.5mm
    \lint_{B_{2\sqrt{t-s}} (x)} \mathbf{e}_\lambda
     \, dy \right)^{1/2}
\notag\\
%4
&
\; \le \; \frac C{R_0} |t-s|^{(d+1)/2}.
\label{INEQ:C-T}
\end{align}
By plying \eqref{INEQ:C-T} and \eqref{INEQ:PR}, we calculate
 the term of $| u_\lambda (t,x) \, - \, u_\lambda (s,x)|$:
\begin{align}
%1
&
| u_\lambda (t,x) \, - \, u_\lambda (s,x) |
\notag\\
%2
&
\; \le \;
 \biggl|  u_\lambda (t,x) \, - \, 
  \avint{B_{\sqrt{t-s}} (x)}  u_\lambda (t,y) \, dy \biggr|
\, + \,
 \biggl|  u_\lambda (s,x) \, - \,
  \avint{B_{\sqrt{t-s}} (x)}  u_\lambda (s,y) \, dy \biggr|
\notag\\
%3
&
\, + \, 
 \avint{B_{\sqrt{t-s}} (x)}  
  | u_\lambda (t,y) \, - \, u_\lambda (s,y) | \, dy
\notag\\
%4
&
\; \le \; \avint {B_{\sqrt{t-s}} (x)} dy
 \lint_0^{|y-x|}
  \left(
   \left| \frac {\partial u_\lambda}{\partial \rho}
          \Bigl(t,x+\rho \frac{y-x}{|y-x|}\Bigr)
   \right|
    \, + \, 
     \left| \frac {\partial u_\lambda}{\partial \rho}
          \Bigl(s,x+\rho \frac{y-x}{|y-x|}\Bigr)
     \right|
      \right) \, d\rho
\notag\\
%5
&
\, + \, \frac C{R_0} |t-s|^{1/2}
\; \le \; \frac {3C}{R_0} |t-s|^{1/2}.
\qquad\qquad\qquad\qquad\qquad\qquad\qquad\qed
\end{align}
%%%%%%%%%%%%%%%%%%%%%%%%%%%%%%%%%%%%%%%%%%%%%%%%%%%%%%%%%%
\begin{Thm}{\rm{(Singular Set).}}\label{THM:Singular}
The set of $\mathbf{sing}$ is a relatively closed set and 
\begin{equation}
\mathcal{H}^{(d)} (\mathbf{sing}) \, = \, 0
\label{EQ:Hausdorff-Estimate}
\end{equation}
holds with respect to the parabolic metric.
\end{Thm}
\noindent{\underbar{Proof of Theorem \ref{THM:Singular}.}$\;$} 
\rm\enspace
\vskip 6pt
$\textbf{sing}$ is a relatively closed set. 
 Indeed, if $z_0$ $\, \in \,$ 
  $\overline{\textbf{sing}} \, \cap \, Q(T)$,
   some sequence $\{ z_\mu \}$ $(\mu = 1,2,\ldots)$
    $\, \subset \, \textbf{sing} \, \cap \, Q(T)$
     satisfies $z_\mu \to z_0$ as $\mu \nearrow \infty$.
More precisely for any positive number $\delta$, 
 there exists a positive integer $\mu_\delta$ 
  such that %$\mu \ge \mu_\delta$ implies 
   $\mathrm{d} (z_\mu, z_0)$ $<$ $\delta$ 
    holds for an arbitrary positive integer $\mu \; \ge \; \mu_\delta$.
From definition on $\textbf{sing}$, 
 for any $R \, > \, \delta$ and any points 
  $z_\mu$ 
   $(\mu = \mu_\delta, \mu_\delta + 1,\ldots)$
    $\in$ $textbf{sing}$ $\cap$ $Q(T)$,
     we obtain
\allowdisplaybreaks\begin{align}
%1
\epsilon_0 & \, \le \; 
 \frac 1 {(R-\delta)^{d}}
  \displaystyle\limsup_{\lambda (\nu) \nearrow \infty}
   \lint_{P_{R-\delta} (z_\mu)} %\kern -5mm
    \mathbf{e}_{\lambda (\nu)} (z) \, dz
\notag
\\
%2
&
\; \le \; 
 \frac 1 {(R-\delta)^{d}}
  \displaystyle\limsup_{\lambda (\nu) \nearrow \infty}
   \lint_{P_R (z_0)} \mathbf{e}_{\lambda (\nu)} (z) \, dz.
\end{align}
By the arbitrariness of $\delta$, 
 passing to the limit of $\delta \searrow 0$, we can say
  $\overline{\textbf{sing}} \cap Q(T)$ $\, \subset \,$ 
   $\mathbf{sing} \cap Q(T)$,
    which provides us with our first assertion.
\par
Next we estimate the size of $\mathbf{sing}$ in the $d$-dimensional
 Hausdorff measure with respect to the parabolic metric. 
Fix a positive integer $n$ and a positive number $R$:
 Set a compact set $Q_{n}$ $=$ $[1/n^2,T-1/n^2]$ $\times$ $\overline{B_{1-1/n} (0)}$ 
  and let $\{ P_{2R_k} (z_k)\}$ $\; \bigl(0 \, < \, 2R_k \, < \, \min
   (R,1/(2n)) \bigr)$ 
     be a cover of $\mathbf{sing}$.
The parabolic version of Vitali covering theorem shows that 
 there is a disjoint finite subfamily $\{ P_{2R_j} (z_j) \}$ $(j \, \in \, \mathcal{J})$ with
\allowdisplaybreaks\begin{align*}
%1
&
\mathbf{sing} \cap Q_{n}
 \, \subset \, \bigcup_{ j \in \mathcal{J}} P_{20R_j} (z_j),
\quad \epsilon_0 R_j^d \; \le \; %\kern -2mm
 \limsup_{\lambda (\nu) \nearrow \infty}
  \lint_{P_{R_j} (z_j)} \mathbf{e}_{\lambda (\nu)} (z) \, dz.
\end{align*}
\par
From Corollary \ref{COR:RPI}, we infer
\allowdisplaybreaks\begin{align}
%1
\epsilon_0 R_j^d & \; \le \;
 \limsup_{\lambda (\nu) \nearrow \infty}
  \lint_{P_{R_j} (z_j)}
   \mathbf{e}_{\lambda (\nu)} (z) \, dz
\notag
\\
%2
&
\le \; \frac C {R_j^2}
 \lint_{P_{2R_j} (z_j)}
  | u (z) \, - \, u_{B_{2R_j} (x_j)} (t) |^2 \, dz
\le \; C
 \lint_{P_{2R_j} (z_j)} | \nabla u (z) |^2  \, dz.
\notag
\end{align}
Thus we obtain
\allowdisplaybreaks\begin{align}
%1
&
\sum_{j=1}^J (20R_j)^d \; \le \; C 
 \lint_{\cup_{j=1}^J P_{2R_j} (z_j)} | \nabla u (z) |^2  \, dz,
\label{INEQ:Absolute-Conti}
\\
%2
&
\mathrm{The \enspace relation \enspace of} \;
 \sum_{j=1}^J (2R_j)^{d+2} 
\; \le \; C R^2 
 \lint_{Q(T)} | \nabla u (z) |^2 \, dz,
\notag
\end{align}
from \eqref{INEQ:Absolute-Conti} and the absolute continuity of the Lebesgue integration, concludes
\begin{equation}
 \mathcal{H}^{(d)} ( \mathbf{sing} \cap Q_n )
  \; \le \; C \lim_{R \searrow 0} \sum_{j=1}^J (20R_j)^{d}
   \; = \; 0.
\end{equation}
\par
By $\lim_{n \to \infty} \mathcal{H}^{(d)}$ 
 $( \textbf{sing} \cap Q_n )$
  $\, = \,$ $\mathcal{H}^{(d)} (\textbf{sing})$,
   we can deduce our assertion. $\qed$
\subsection{\enspace Compactness}
\par
We have seen in the previous section that
 the GLHF $u_\lambda$ converges weakly to a WHHF $u$
  (modulo a subsequence of $\lambda$).
   Theorem 3 in  L.~C.~Evans \cite[p.39]{evans-88} expounds
    that the differentials of the GLHF, $\nabla u_\lambda$ does strongly
     to them of the WHHF, $\nabla u$ in $L^p (Q(T) ; \mathbb{R}^{d(D+1)})$ with $1 < p < 2$.
But it doesn't suffice to prove that
 the WHHF $u$ satisfies the monotonicity for the scaled energy,
  Corollary \ref{COR:Mon}.
Thus we demonstrate the strong convergencity
 of $\{ u_{\lambda (\nu)} \}$ $( \nu = 1,2, \ldots)$ to a WHHF $u$ 
  in $H_{\mathrm{loc}}^{1,2}$-topology as $\lambda (\nu) \nearrow \infty$.
   The estimates of $\mathcal{H}^{(d)} (\mathbf{sing})$ $\, = \,$ $0$ in Theorem
    \ref{THM:Singular} plays a crucial role in the proof.
\begin{Thm}{\rm{(Strong Convergencity of Gradients of \textrm{WHHF}).}}
\label{THM:Strong-Convergence-Gradient}
For a suitable subsequence of $\{ \lambda (\nu) \}$
 still denoted by $\{ \lambda (\nu) \}$
  $( \nu = 1,2, \ldots )$, a sequence of the gradients of the GLHF,
   $\{ \nabla u_{\lambda (\nu)}\}$ $( \nu = 1,2, \ldots)$ converges
    strongly to the gradients of the WHHF in $L_{\loc}^2 (Q(T) ; \mathbb{R}^{d(D+1)})$.
\end{Thm}
\noindent{\underbar{Proof of Theorem \ref{THM:Strong-Convergence-Gradient}.}} 
\rm\enspace
\vskip 6pt
Set $Q_n$ $=$ $[1/n^2, T-1/n^2]$ $\times$ $\overline{{B}_{1-1/n} (0)}$ 
 $( n = 1,2,\ldots )$ and
  fix any compact sets
   $Q_n$ $\subset\!\subset$ $Q_{2n}$ $\subset\!\subset$ $Q (T)$ and any positive integer $k$.
By means of $\mathcal{H}^{(d)} (\mathbf{sing})$ $\, = \,$ $0$,
 we first see that we can choose up 
  two finite sets of cylinders $\{ P_{r_{i,k}} (z_{i,k}) \}$
   and $\{ P_{\rho_{j,k}} (z'_{j,k}) \}$
    $( i = 1,2,\ldots, I_k \, ; \, j = 1,2,\ldots, J_k )$ satisfy that
\allowdisplaybreaks\begin{align}
%1
&
\mathbf{sing} \cap Q_{2n}
 \, \subset \, \bigcup_{i=1}^{I_k} P_{r_{i,k}} (z_{i,k}), \;
\quad \sum_{i=1}^{I_k} r_{i,k}^d
 \; \le \; \frac 1k,
\label{INEQ:claim-1}
\\
%2
&
Q_{2n} \setminus \bigcup_{i=1}^{I_k} P_{r_{i,k}} (z_{i,k})
 \, \subset \, 
  \bigcup_{j=1}^{J_k} P_{\rho_{j,k}} (z'_{j,k})
   \, \subset\!\subset \, Q(T)
\label{INEQ:claim-2}
\\
%3
&
\mathrm{with} \quad
 \overline{\mathcal{M}} ( P_{g(\rho_{j,k})} (z'_{j,k}))
  \; < \; \epsilon_0
\label{INEQ:claim-3}
\end{align}
where the function $g$ is the positive function appeared in Theorem \ref{THM:PR}.
In addition by a diagonal argument, we find that we can pick up a subsequence
 $\{ \lambda ( \nu ( l ) ) \}$ $( l = 1,2,\ldots)$
  of $\{ \lambda ( \nu ) \}$ $( \nu = 1,2,\ldots)$
   with the following properties:
\renewcommand{\labelenumi}{(\roman{enumi})}
\begin{enumerate}
\item
$\lambda ( \nu ( l ) )$ $\ge$ $l$,
\item
for any positive integer $k$,
 $| u_{\lambda (\nu(l)) } (z) \, - \, u (z) |$ $\; \le \;$ $1/k$
  holds on any point $z$ $\in$ 
   $Q_{2n} \setminus \bigcup_{i=1}^{I_k} P_{r_{i,k}} (z_{i,k})$
    whenever $l$ $\ge$ $k$.
\end{enumerate}
Indeed from \eqref{INEQ:claim-3}, it follows that
 for some $\lambda (\nu_{j,k})$ depending on $j$ and $k$,
\begin{equation}
\frac 1 {g(2\rho_{j,k})^d}
 \lint_{P_{g(2\rho_{j,k})} (z'_{j,k})}
  \mathbf{e}_{\lambda{(\nu)}} (z) \, dz
   \; < \; \epsilon_0
\notag
\end{equation}
holds if $\lambda (\nu)$ is more than or equal to 
 $\lambda (\nu_{j,k})$. 
Then Theorem \ref{THM:PR} and  Lemma \ref{LEM:Cont-Time} can read
\begin{equation}
| u_{\lambda (\nu)} (z_1) \, - \, u_{\lambda (\nu)} (z_2) |
 \; \le \; \frac {C J_k}{\rho_{j,k}} d(z_1,z_2)
\label{INEQ:Strong-3}
\end{equation}
for any points $z_1$ and $z_2$ $\in$ 
 $Q_{2n} \setminus \cup_{i=1}^{I_k} P_{\rho_{i,k}} (z_{i,k})$.
\par
When we set $\bar{\nu}_k$ $=$ $\max_{j=1,2,\ldots,J_k}$
 $\nu_{j,k}$ and $\underline{\rho}_k$ $=$ $\min_{j=1,2,\ldots,J_k}$
  $\rho_{j,k}$ in \eqref{INEQ:Strong-3},
Ascoli-Arzela's theorem claims that 
 there exists a subsequence $\{ \lambda (\nu(l))\}$
  $(l=1,2,\ldots)$ of $\{ \lambda (\nu)\}$
   $(\nu = 1,2,\ldots)$ such that
    for some $\nu(k)$ more than or equal to 
     $\bar \nu_k$ and $k$
      if $\nu (l)$ $\ge$ $\nu (k)$,
\begin{equation}
| u_{\lambda (\nu (l) ) } (z) \, - \, u (z) | \; \le \; \frac 1k
\label{INEQ:Strong-6}
\end{equation}
holds on any point $z$ $\in$ $Q_{2n} \setminus \bigcup_{i=1}^{I_k} P_{r_{i,k}} (z_{i,k})$,
 where the mapping $u$ is the WHHF constructed in Theorem
  \ref{THM:Existence-HF}.
\par
Next we do exactly the same  procedure above for $k+1$ instead of $k$.
 So we can select a number $\lambda (\nu (k+1))$ 
  from $\{ \lambda (\nu(l))\}$ $(l=1,2,\ldots)$
   more than or equal to $\lambda (\nu (k))$ and $k+1$ respectively satisfying

\allowdisplaybreaks\begin{align}
%1
&
| u_{\lambda (\nu ( k+1)) } (z) \, - \, u (z) | \; \le \; \frac
 {1}{k+1},
\notag
\\
%2
&
| u_{\lambda (\nu ( k+1)) } (z) \, - \, u (z) | \; \le \; \frac
 {1}{k},
\label{INEQ:Strong-7}
\end{align}
holds on any point $z$ $\in$
 $Q_{2n} \setminus \bigcup_{i=1}^{I_{k+1}} P_{r_{i,k+1}} (z_{i,k+1})$
  and
   $Q_{2n} \setminus \bigcup_{i=1}^{I_{k}} P_{r_{i,k}} (z_{i,k})$.
\par
%Adopt $\lambda (\nu (k) )$ as $\lambda (k)$;
 By an induction, we can choose a subsequence $\{ \lambda (\nu(k)) \}$
 $(k=1,2,\ldots)$ of $\{ \lambda \}$ $( \lambda > 0)$  satisfying
  monotone nondecreasing with respect to $k$, and
\allowdisplaybreaks\begin{align}
%1
&
\lambda (\nu(l)) \; \ge \; \max ( \lambda (\nu(k)), l),
\notag
\\
%2
&
| u_{\lambda (\nu(l))} (z) \, - \, u (z) | \; \le \; \frac 1k 
 \qquad {\mathrm{on}} \quad 
  Q_{2n} \setminus \bigcup_{i=1}^{I_{k}} P_{R_{i,k}} (z_{i,k})
\notag
\end{align}
holds for any integers $k$ and $l$ with $k \le l$.
\par
An adaptation of $\lambda (\nu(k))$ to $\lambda (k)$
 shows our claim.
\par
We shall confirm that the sequence of $\{ \nabla u_{\lambda (k)} \}$ $( k = 1,2, \ldots, )$
 does converge {\it{strongly}} to the gradient of the WHHF $u$
  appeared in Theorem \ref{THM:Existence-HF} in
   $L_{\loc}^2 (Q_{2n} ; \mathbb{R}^{d(D+1)})$ as $k \nearrow \infty$.
    To this end, fix any positive integer $k$ and let
     $\lambda (l)$ be more than or equal to $\lambda (k)$.
Take the difference between \eqref{EQ:GL} in $\lambda$ $=$ $\lambda (l)$
 and \eqref{EQ:HHF},
\allowdisplaybreaks\begin{align}
%1
&
\lint_{Q_{2n}}
 \Bigl\la \frac{\partial}{\partial t}
  ( u_{\lambda (l)} \, - \, u ), \phi 
    \Bigr\ra \, dz
\, + \,
 \lint_{Q_{2n}}
  \la \nabla ( u_{\lambda (l)} \, - \, u ),
   \nabla \phi \ra \, dz
\notag
\\
%2
&
\; = \; - \lint_{Q_{2n}}
 | \nabla u|^2 \la u,\phi \ra \, dz
\, + \,
 \lint_{Q_{2n}} 
  \lambda (l)^{1-\kappa }
   \bigl( 1 \, - \, | u_{\lambda (l)} |^2 \bigr) 
    \la u_{\lambda (l)}, \phi \ra \, dz
\label{EQ:Diff-EQ}
\end{align}
for a map $\phi \in C_0^\infty (Q_{2n} ; \mathbb{R}^{D+1})$.
Choose a smooth function $\eta_1$ satisfying
\begin{equation}
\eta_1 (z) \; = \;
\begin{cases}
%1.1
1 & \mathrm{on} \; Q_{n},
\\
%1.2
0 & \mathrm{off} \; Q_{2n}
\end{cases}
\end{equation}
with $0 \,\le\, \eta_1 \,\le\, 1$, $| \nabla \eta_1 |$ $\le$ $2n$,
 $| \triangle \eta_1 |$ $+$
  $|\partial \eta_1 / \partial t|$ $\le$ $16n^2$.
\par
We substitute a smooth approximation of 
 $( u_{\lambda (l)} \, - \, u ) \eta_1$ for $\phi$
  in \eqref{EQ:Diff-EQ}.
   After passing to the limit, we obtain
\allowdisplaybreaks\begin{align}
%1
&
\lint_{Q_{2n}}
 | \nabla ( u_{\lambda (l)} \, - \, u )|^2 
  \eta_1 \, dz
\; =\; \frac 12 \lint_{Q_{2n}}
 | u_{\lambda (l)} \, - \, u |^2 \,
  \Bigl( \frac{\partial \eta_1}{\partial t} 
   \, + \, 
    \triangle \eta_1 \Bigr)
     \, dz
\notag
\\
%3
&
\, - \, \lint_{Q_{2n}}
 |\nabla u |^2 | 
  \la u, u_{\lambda (l)} \, - \, u \ra \eta_1 \, dz
\notag
\\
%4
&
\, + \, \lint_{Q_{2n}}
 \lambda (l)^{1-\kappa }
  ( 1 \, - \, |u_{\lambda (l)}|^2 )
    \la u_{\lambda (l)}, u_{\lambda (l)} \, - \, u \ra \eta_1 \, dz.
\label{INEQ:4_0}
\end{align} 
\par
From now on, we compute to pass to the limit 
 $\lambda (l)$ $\nearrow$ $\infty$
  on the right-hand side in \eqref{INEQ:4_0}.
By using the strong convergencity of $u_{\lambda (l)}$ in 
 $L^2 (Q(T) ; \mathbb{R}^{{D+1}})$, i.e. Theorem \ref{THM:Convergence},
  we can calculate the first term as follows:
\allowdisplaybreaks\begin{align}
%1
&
\limsup_{\lambda (l) \nearrow \infty} 
 \lint_{Q_{2n}}
  | u_{\lambda (l)} \, - \, u |^2 \,
   \Bigl( \frac{\partial \eta_1}{\partial t}
    \, + \,
     \triangle \eta_1 \Bigr)
      \, dz
\; = \; 0.
\label{INEQ:GRA_1}
\end{align}  
\par
Next we estimate the second term on the right-hand side:
 noting $|u_{\lambda (l)}|$ $\, \le \, 1$, 
  thanks to the dominated Lebesgue convergence theorem
   and Theorem \ref{THM:Convergence}, we obtain
\begin{equation}
\limsup_{\lambda (l) \nearrow \infty}
 \lint_{Q_{2n}}
  | \nabla u |^2 | u_{\lambda (l)} \, - \, u | \, dz
   \; = \; 0.
\label{INEQ:GRA_2}
\end{equation}
Finally we asses the third term. We decompose it into
\allowdisplaybreaks\begin{align}
%1
&
\lint_{Q_{2n}}
 {\lambda (l)^{1-\kappa }}
  ( 1 \, - \, | u_{\lambda (l)} |^2 )
   \la u_{\lambda (l)}, u_{\lambda (l)} \, - \, u \ra \eta_1 \, dz
\notag
\\
%2
&
\; = \;
 \lint_{\cup_{i=1}^{I_k} P_{r_{i,k}} (z_{i,k})} % \kern-3mm
  {\lambda (l)^{1-\kappa }}
  ( 1 \, - \, | u_{\lambda (l)}|^2) 
    \la u_{\lambda (l)}, u_{\lambda (l)} \, - \, u \ra \eta_1 \, dz
\notag
\\
%3
&
\, + \, 
 \lint_{Q_{2n}\setminus\cup_{i=1}^{I_k} P_{r_{i,k}} (z_{i,k})}
   \lambda (l)^{1-\kappa }
   ( 1 \, - \, | u_{\lambda (l)}|^2) 
     \la u_{\lambda (l)}, u_{\lambda (l)} \, - \, u \ra \eta_1 \, dz.
\label{INEQ:4_1}
\end{align}
We majorize the first term in \eqref{INEQ:4_1} as follows:
 Recall a way of choosing a sequence of cylinder $\{ P_{r_{i,k}} (z_{i,k}) \}$
  $( i = 1,2,\ldots, I_k )$
   with \eqref{INEQ:claim-1}
    for any positive integer $k$
and set a certain smooth cut off functions $\phi_{i,k}$ with
 $0 \le \phi_{i,k} \le 1$, $| \nabla \phi_{i,k} |$ $\le$ $2/r_{i,k}$,
  $| \triangle \phi_{i,k} |$ $\le$ $4/r_{i,k}^2$ and
   $| \partial \phi_{i,k} / \partial t| $ $\le$ $2/r_{i,k}^2$
    with
\begin{equation*}
\phi_{i,k} \; = \;
\begin{cases}
%1.1
1 & \; \mathrm{in} \quad P_{r_{i,k}} (z_{i,k}),
\\
%1.2
0 & \; \mathrm{outside} \quad P_{2r_{i,k}} (z_{i,k}).
\end{cases}
\end{equation*}
Multiply \eqref{EQ:GL} in $\lambda = \lambda (l)$ by $u_{\lambda (l)}$ $\phi_{i,k}$
 and integrate it on $P_{r_{i,k}} (z_{i,k})$ to observe
\allowdisplaybreaks\begin{align}
%1
&
\lint_{P_{r_{i,k}} (z_{i,k})}
 \lambda (l)^{1-\kappa }
  ( 1 \, - \, | u_{\lambda (l)} |^2 ) \, dz
\notag
\\
%2
&
\; \le \; \frac C{r_{i,k}^2} 
 \lint_{P_{2r_{i,k}} (z_{i,k})}
  ( 1 \, - \, | u_{\lambda (l)} |^2 ) \, dz
\, + \, C
 \lint_{P_{2r_{i,k}} (z_{i,k})}
  \mathbf{e}_{\lambda (l)} \, dz.
\label{INEQ:4_2}
\end{align}
Set $z_{i,k}$ $=$ $(t_{i,k},x_{i,k})$ and
 recall $d_n$ $=$ $\mathrm{d} (Q_{n},\partial Q(T))$$ ( = 1/n )$;
  Thus using Theorem \ref{THM:Energy-Estimate} and Corollary \ref{COR:Mon},
   the second term in \eqref{INEQ:4_2} can be evaluated
\allowdisplaybreaks\begin{align}
%1
&
\lint_{P_{r_{i,k}} (z_{i,k})}
 \mathbf{e}_{\lambda (l)} \, dz
\; = \;
 \lint_{t_{i,k} + {20 r_{i,k}^2}/3 - {32 r_{i,k}^2}/3}
      ^{t_{i,k} + {20 r_{i,k}^2}/3 - {8 r_{i,k}^2}/3}
  dt
   \lint_{B_{r_{i,k}} (x_{i,k})}
     \mathbf{e}_{\lambda (l)} \, dx
\notag\\
%2
&
\; \le \; \frac{C r_{i,k}^d}{r_{i,k}^d}
\lint_{t_{i,k} + {20 r_{i,k}^2}/3 - {32 r_{i,k}^2}/3}
      ^{t_{i,k} + {20 r_{i,k}^2}/3 - {8 r_{i,k}^2}/3}
  dt
   \lint_{B_{r_{i,k}} (x_{i,k})}
    \mathbf{e}_{\lambda (l)} 
     \, \exp \biggl( 
      \frac {| x - x_{i,k} |^2}{4(t - (t_{i,k} + 20r_{i,k}^2/3 ))} \biggr)
       \, dx
\notag\\
%3
&
\; \le \; C r_{i,k}^d
 \biggl(
  \frac 1{d_n^d}
   \lint_{t_{i,k} + 20 r_{i,k}^2/3 - d_n^2}
        ^{t_{i,k} + 20 r_{i,k}^2/3 - (d_n/2)^2} \, dt
     \lint_{\mathbb{B}^d} \mathbf{e}_{\lambda (l)}
      \, dx
\, + \, C r_{i,k}^2 \biggr)
\notag\\
%4
&
\; \le \; C r_{i,k}^d.
\label{INEQ:4_3}
\end{align}
Then from \eqref{INEQ:4_1}, \eqref{INEQ:4_2} and \eqref{INEQ:4_3}, 
 we obtain
\allowdisplaybreaks\begin{align}
%1
&
\lint_{\cup_{i=1}^{I_k} P_{r_{i,k}} (z_{i,k})} %\kern-3mm
 {\lambda (l)^{1-\kappa }}
  ( 1 \, - \, | u_{\lambda (l)}|^2) 
    \la u, u_{\lambda (l)} \, - \, u \ra \eta_1 \, dz
\notag
\\
%2
&
\; \le \; C \sum_{i=1}^{I_k}
 \lint_{P_{r_{i,k}} (z_{i,k})}
  \mathbf{e}_\lambda (z) \, dz
\, + \, \sum_{i=1}^{I_k}
 \frac C {r_{i,k}^2}
  \lint_{P_{r_{i,k}} (z_{i,k})} 
   ( 1 \, - \, | u_{\lambda (l)}|^2) \, dz
\notag\\
%3
&
\; \le \; C \sum_{i=1}^{I_k}  r_{i,k}^d
 \; \le \; \frac Ck.
\label{INEQ:4_4}
\end{align}
On the other hand, since $|u_{\lambda (l)} \, - \, u |$ $\le$  $1/k$ on
 $Q_{2n}\setminus\cup_{i=1}^{I_k} P_{r_{i,k}} (z_{i,k})$,
  likewise \eqref{INEQ:4_2},
   the second term of \eqref{INEQ:4_1} becomes
\allowdisplaybreaks\begin{align}
%1
&
\lint_{Q_{2n} \setminus \cup_{i=1}^{I_k} P_{r_{i,k}} (z_{i,k})}
 \lambda (l)^{1-\kappa }
  ( 1 \, - \, | u_{\lambda (l)} |^2 ) 
   \la u_{\lambda (l)}, u_{\lambda (l)} \, - \, u \ra \eta_1 \, dz
\notag
\\
%2
&
\; \le \; \frac Ck
 \lint_{Q_{2n} \setminus \cup_{i=1}^{I_k} P_{r_{i,k}} (z_{i,k})}
  \lambda (l)^{1-\kappa }
   ( 1 \, - \, | u_{\lambda (l)} |^2 ) \eta_1 \, dz
\notag
\\
%3
&
\; \le \; \frac Ck
 \lint_{Q_{2n}}
  \lambda (l)^{1-\kappa }
   ( 1 \, - \, | u_{\lambda (l)} |^2 ) \eta_1 \, dz
\notag\\
%4
&
\; \le \; \frac Ck
 \biggl(
  \frac C{d_n^2}
   \lint_{Q(T)} ( 1 \, - \, | u_{\lambda (l)} |^2 ) \, dz
    \, + \, T 
     \lint_{\mathbb{B}^d} | \nabla u_0 |^2 \, dx
      \biggr)
\label{INEQ:4_5}
\end{align}
as long as $\lambda (l)$ $\ge$ $\lambda (k)$.
Taking supremum limit in \eqref{INEQ:4_4} and \eqref{INEQ:4_5}, we deduce 
\begin{equation}
\limsup_{\lambda (l) \nearrow \infty}
 \lint_{Q_{2n}}
  \lambda (l)^{1-\kappa }
   ( 1 \, - \, | u_{\lambda (l)} |^2 ) 
    \la u_{\lambda (l)},
      u_{\lambda (l)} \, - \, u
       \ra \eta_1 \, dz
\; \le \; \frac Ck.
\label{INEQ:GRA_3}
\end{equation}
Since $k$ is an arbitrary positive integer,
 our conclusion follows from \eqref{INEQ:4_0}, \eqref{INEQ:GRA_1}, \eqref{INEQ:GRA_2}
  and \eqref{INEQ:GRA_3}.
   \qed

%% file: main_hhf_arch.tex
%
%#! platex hhf
%
\setcounter{chapternumber}{4}\setcounter{equation}{0}
\renewcommand{\theequation}%
           {\thechapternumber.\arabic{equation}}
\section{\enspace Proof of Main Theorems}
\par
By making the best of a few ingredients 
 and properties on the WHHF and the GLHF,
  this chapter establishes Theorem \ref{THM:Main-1} and Theorem \ref{THM:Main-2} in Chapter \ref{SEC:Intro}.
\subsection{Proof of Theorem \ref{THM:Main-1}}
Let $\{ \lambda (k) \}$ $( k \, = \, 1,2,\ldots )$
 be a subsequence of $\lambda$ be chosen in Theorem
  \ref{THM:Strong-Convergence-Gradient}.
Theorem \ref{THM:Existence-HF} tells us that
 the subsequence of the GLHF $\{ u_{\lambda (k)}\}$
  $( k = 1, 2, \ldots )$ converges to a WHHF $u$
   in $L^2 (Q(T) ; \mathbb{R}^{D+1} )$.
\par
Next we discuss a partial regularity on the WHHF
 constructed above:
On account of Theorem \ref{THM:Energy-Estimate} and
 Theorem \ref{THM:Strong-Convergence-Gradient},
  we obtain
\allowdisplaybreaks\begin{align}
%1
&
\mathbf{sing} \, = \, \bigcap_{R > 0}
 \bigl\{ z_0 \in Q(T) \, ; \, 
  \frac 1 {2R^d} \kern-2.5pt \lint_{P_R (z_0)} \kern-2.5pt
   | \nabla u |^2 \, dz
    \, \ge \, \epsilon_0, \,
P_R (z_0) \, \subset\!\subset \, Q(T)
\bigr\}
\label{EQ:Sing-2}
\\
%2
&
\mathbf{reg} \; = \; Q(T) \setminus \mathbf{sing},
\label{EQ:Reg-2}
\end{align}
where a number $\epsilon_0$ is a positive constant appeared in Theorem
\ref{THM:PR}.
From Theorem \ref{THM:Singular}, we see that \textbf{sing}
 is relatively closed. 
First of all, we measure the size of \textbf{sing}.
 Recall that by Lemma \ref{LEM:RH}, 
  the differentials $\nabla u$ of the WHHF belongs to $L_{\loc}^{2q_0}$ $(Q(T) ;
   \mathbb{R}^{d(D+1)} )$
    for some positive number $q_0$ greater than $1$.
So accordingly the inclusion
\begin{equation*}
\mathbf{sing} \; \subset \;
 \bigcap_{R > 0}
  \Bigl\{ z_0 \in Q(T) \, ; \, 
   \frac 1{R^{d-2(q_0-1)}} 
    \kern-2.5pt \lint_{P_R (z_0)} \kern-2.5pt
    | \nabla u |^{2q_0} \, dz \,
     \, \ge \, ( C \epsilon_0)^{q_0}, \, 
      P_R (z_0) \, \subset\!\subset \, Q(T)
\Bigr\}
\end{equation*}
enjoys
\begin{equation}
\mathcal{H}^{d-2(q_0-1)} (\mathbf{sing})
 \; < \; \infty.
\label{INEQ:Sing}
\end{equation}
Next we show that if $z_0$ $\in$ $\textbf{reg}$,
 the WHHF $u$ is smooth on some parabolic cylinder $P_{R_0/4} (z_0)$
  by plying Lady\v{z}henskaya,~O.~A., Solonnikov,~V.~A.,
   Ural'ceva,~N.~N, \cite{ladyzhenskaya-solonnikov-uralceva}.
Invoke Theorem \ref{THM:PR} and Lemma \ref{LEM:Cont-Time}
 to arrive at
\allowdisplaybreaks\begin{align}
%1
&
| \nabla u(z) | \; \le \; \frac C{R_0}
 \quad \mathrm{and} \quad
  | u(t,x) \, - \, u(s,x) | \; \le \; 
   \frac {C}{R_0} |t-s|^{1/2}
\end{align}
on any point $z$, $(t,x)$ and $(s,x)$ $\in$ $P_{R_0/2} (z_0)$
 for some positive number $R_0$.
We prepare a smooth cut-off function $\eta_{R_0}$ given by
\begin{equation}
\eta_{R_0} (z) \; = \;
\begin{cases}
%1.1
1 & \mathrm{in} \quad P_{R_0/4} (z_0),
\\
%1.2
0 & \mathrm{off} \quad P_{R_0/2} (z_0)
\end{cases}
\end{equation}
with $0 \le \eta_{R_0} \le 1$.
Then the mapping $u \eta_{R_0}$ satisfies the following system:
\begin{align}\allowdisplaybreaks
%1
&
\Bigl( 
 \frac {\partial}{\partial t} \, - \, \triangle
  \Bigr)
   (u \eta_{R_0})
    \; = \; | \nabla u |^2 u \eta_{R_0}
\notag\\
%2
&
\, - \, 2 \la \nabla \eta_{R_0}, \nabla u \ra
 \, - \, u \triangle \eta_{R_0}
  \quad \mathrm{in} \quad P_{R_0/2} (z_0).
\label{EQ:Modified-HHF}
\end{align}
By applying Lady\v{z}henskaya,~O.~A., Solonnikov,~V.~A., Ural'ceva,~N.~N.
 \cite[p.341, Theorem 9.1]{ladyzhenskaya-solonnikov-uralceva},
  we know that the mapping of $u \eta_{R_0}$ belongs to $H^{2,q}$ $(P_{R_0/2} (z_0))$
   for any number $q$ greater than $1$ because
    $| \nabla u |^2 u \eta_{R_0}$ 
     $-$ $2 \la \nabla \eta_{R_0}, \nabla u \ra$
      $-$ $u \triangle \eta_{R_0}$
       is bounded on $P_{R_0/2} (z_0)$.
Then employ \cite[p.80, Lemma 3.3]{ladyzhenskaya-solonnikov-uralceva} to verify
 $\nabla u$ $\in$ $C^{\alpha_0}$ $(P_{R_0/4} (z_0))$
   where $\alpha$ is any arbitrary positive number less than $1$.
So by having Schauder estimate at our disposal, 
 it is shown $u$ $\in$ $C^{2+\alpha}$ $(P_{R_0/4} (z_0))$.
  For Schauder estimate, we refer to 
   \cite[p.320, Theorem 5.2]{ladyzhenskaya-solonnikov-uralceva}.
In the light of a boot strap argument,
  we finally get the smoothness of $u$ on $P_{R_0/4} (z_0)$.
\par
We readily see that 
 the WHHF $u$ satisfies a global energy inequality (i),
  a monotonicity for the scaled energy (ii) and
   a reverse Poincar\'{e} inequality (iii).
In fact, the first and third inequalities are established 
 by using Theorem \ref{THM:Convergence} 
  and taking the limit inferior
   of $\lambda (k)$ $\nearrow$ $\infty$ in \eqref{INEQ:Energy-Estimate1}
    in Theorem \ref{THM:Energy-Estimate} and \eqref{INEQ:RPI} in
     Corollary \ref{COR:RPI}.
While the second is proved as follows:
 For a point $z_0$ $=$ $(t_0,x_0)$,
  set $d_0$ $=$ $\dist (x_0, \partial \mathbb{B}^d)$.
Recall Corollary \ref{COR:Mon} and note
\allowdisplaybreaks\begin{align}
&
\frac 1 {R_2^d}
 \lint_{t_0-(2R_2)^2}^{t_0-R_2^2} \, dt
  \lint_{\mathbb{B}^d}
   \mathbf{e}_{\lambda (k)}
    \exp \Bigl( \frac {|x-x_0|^2}{4(t-t_0)}\Bigr)
     \, dx
\notag
\\
%2
&
\; \le \; \frac 1 {R_2^d}
 \lint_{t_0-(2R_2)^2}^{t_0-R_2^2} \, dt
  \lint_{{B}_{1-d_0} (0)}
   \mathbf{e}_{\lambda (k)}
    \exp \Bigl( \frac {|x-x_0|^2}{4(t-t_0)}\Bigr)
     \, dx
\notag
\\
%2
&
\, + \,  \frac 1 {R_2^d}
 \lint_{t_0-(2R_2)^2}^{t_0-R_2^2} \, dt
  \lint_{\mathbb{B}^d \setminus {B}_{1-d_0} (0)}
   \mathbf{e}_{\lambda (k)}
    \exp \Bigl( \frac {|x-x_0|^2}{4(t-t_0)}\Bigr)
     \, dx.
\label{INEQ:Mon-Right}
\end{align}
\par
We pass to the limit $\lambda (k)$ $\nearrow$ $\infty$
 in above \eqref{INEQ:Mon-Right}:
  By  using Theorem \ref{THM:Strong-Convergence-Gradient}  
   for the first term
    and Theorem \ref{THM:Energy-Estimate}
     for the second term,
we obtain
\allowdisplaybreaks\begin{align}
&
\frac 1 {R_2^d}
 \limsup_{\lambda(k) \nearrow \infty}
  \lint_{t_0-(2R_2)^2}^{t_0-R_2^2} \, dt
   \lint_{\mathbb{B}^d}
    \mathbf{e}_{\lambda (k)}
     \exp \Bigl( \frac {|x-x_0|^2}{4(t-t_0)}\Bigr)
      \, dx
\notag
\\
%2
&
\; \le \; \frac 1 {R_2^d}
 \limsup_{\lambda(k) \nearrow \infty}
  \lint_{t_0-(2R_2)^2}^{t_0-R_2^2} \, dt
   \lint_{{B}_{1-d_0} (0)}
    \mathbf{e}_{\lambda (k)}
     \exp \Bigl( \frac {|x-x_0|^2}{4(t-t_0)}\Bigr)
      \, dx
\notag
\\
%2
&
\, + \,  \frac 1 {R_2^d}
 \limsup_{\lambda(k) \nearrow \infty}
  \lint_{t_0-(2R_2)^2}^{t_0-R_2^2} \, dt
   \lint_{\mathbb{B}^d \setminus {B}_{1-d_0} (0)}
    \mathbf{e}_{\lambda (k)}
     \exp \Bigl( \frac {|x-x_0|^2}{4(t-t_0)}\Bigr)
      \, dx
\notag
\\
%3
&
\; \le \;
 \frac 1 {R_2^d}
  \lint_{t_0-(2R_2)^2}^{t_0-R_2^2} \, dt
   \lint_{{B}_{1-d_0} (0)}
    | \nabla u |^2 
     \exp \Bigl( \frac {|x-x_0|^2}{4(t-t_0)}\Bigr)
      \, dx
\notag
\\
%4
&
\, + \, C(d_0) R_2
 \lint_{\mathbb{B}^d}
  | \nabla u_0 |^2 \, dx,
\label{INEQ:Mon-Right-2}
\end{align}
which assert \eqref{INEQ:Main-Energy-Estimate-2}.
$\qed$
\vskip 6pt
\subsection{Proof of Theorem \ref{THM:Main-2}}
\rm\enspace\noindent
First, we prove that the limiting mapping $u$ is a WHHF and
 it satisfies (\romannumeral 1) on $Q(T)$.
  Set a positive number $h$ sufficiently small;
   From Theorem \ref{THM:Energy-Estimate}, we infer
\allowdisplaybreaks\begin{align}
%1
&
\fint_{T_1-h}^{T_1+h} \, dt_1
 \fint_{T_2-h}^{T_2+h} \, dt_2
  \lint_{t_1}^{t_2} \, dt
   \lint_{\mathbb{B}^d}
    \left| \frac{\partial u}{\partial t } (z) \right|^2 \, dx
\notag
\\
%2
&
\, + \, 
 \fint_{T_2-h}^{T_2+h} \, dt 
  \lint_{\mathbb{B}^d}
   \mathbf{e}_\lambda (t) \, dx
\; \le \; 
 \fint_{T_1-h}^{T_1+h} \, dt
  \lint_{\mathbb{B}^d}
   \mathbf{e}_\lambda (t) \, dx
\label{INEQ:Energy-Estimate-Average}
\end{align}
holds for any positive numbers $T_1$ and $T_2$ 
 with $0$ $<$ $T_1$ $\le$ $T_2$ $<$ $T$.
Our sending $\lambda$ $\nearrow$ $\infty$ of
 \eqref{INEQ:Energy-Estimate-Average},
  in which Theorem \ref{THM:Energy-Estimate}
   and Theorem \ref{THM:Strong-Convergence-Gradient}
    are implemented,
we infer on the same $T_1$ and $T_2$,
\allowdisplaybreaks\begin{align}
%1
&
\lint_{T_1+h}^{T_2-h} \, dt
 \lint_{\mathbb{B}^d}
  \left| \frac{\partial u}{\partial t } (z) \right|^2 \, dx
\notag
\\
%2
&
\, + \, \frac 12
 \fint_{T_2-h}^{T_2+h} \, dt 
  \lint_{\mathbb{B}^d}
   | \nabla u (t) |^2 \, dx
\; \le \; 
 \fint_{T_1-h}^{T_1+h} \, dt
  \lint_{\mathbb{B}^d}
   | \nabla u (t) |^2 \, dx.
\label{INEQ:Energy-Estimate-Average-2}
\end{align}
Thus we possibly pass to the limit 
 $h \searrow 0$ at any Lebesgue points $T_1$ and $T_2$ 
  with $0$ $<$ $T_1$ $\le$ $T_2$ $<$ $T$
   to conclude
\begin{align}
%1
&
\lint_{T_1}^{T_2} \, dt
 \lint_{\mathbb{B}^d}
  \left| \frac{\partial u}{\partial t } (z) \right|^2 \, dx
\, + \, \frac 12
 \lint_{\mathbb{B}^d}
  | \nabla u (T_2) |^2 \, dx
\; \le \; \frac 12
 \lint_{\mathbb{B}^d}
  | \nabla u (T_1) |^2 \, dx.
\label{INEQ:Energy-Estimate-Average-3}
\end{align}
Recalling \eqref{EQ:GL} in Remark \ref{REM:EQ},
 by \eqref{INEQ:Energy-Estimate-Average-3}, Theorem \ref{THM:Convergence} and Theorem \ref{THM:Existence-HF}, 
  we construct a WHHF $u$ $\in$ $V (Q(T);\mathbb{S}^D)$ with
\allowdisplaybreaks\begin{align}
%1
&
\lint_0^T \, dt
 \lint_{\mathbb{B}^d}
  \left| \frac{\partial u}{\partial t } (z) \right|^2 \, dx
\, + \, \frac 12 \,
 \lint_{\mathbb{B}^d}
 | \nabla u (T) |^2 \, dx
\; \le \; 
 \frac 12 \lint_{\mathbb{B}^d} | \nabla u_0 |^2 \, dx,
\notag
\\
%2
&
\lint_{Q(T)}
 \Bigl\langle \frac{\partial u}{\partial t}, \phi
  \Bigr\rangle \, dz
   \, + \, \lint_{Q(T)}
    \langle \nabla u, \nabla \phi \rangle \, dz
      \; = \; \lint_{Q(T)}
      | \nabla u |^2 \langle u, \phi \rangle \, dz
\notag
\end{align}
for any map $\phi$ $\in$ $C^0 ( 0,T \, ; \, C_0^\infty (\mathbb{B}^d ; \mathbb{R}^{D+1}))$.
\par
By the same procedure above on $(T,2T]$, 
 we can extend the WHHF in $Q(T)$ to it in $Q(2T)$
and a repeat argument permits us to comprise a WHHF $u$ $\in$ $V (Q(\infty))$.
\par
Next we prove the constancy property. 
 A combination of Theorem \ref{THM:EDE} with 
  Theorem \ref{THM:Convergence}
   and Theorem \ref{THM:Existence-HF} yields
\begin{equation}
\lint_{\mathbb{B}^d} | \nabla u (t) |^2 \, dx \,
 \; \le \; 2 e^{-(d-2) t} \,
  \lint_{\mathbb{B}^d} | \nabla u_0 |^2 \, dx
   \quad \mathrm{for} \; \mathrm{all} \; t \in (0,\infty).
\label{INEQ:Strong}
\end{equation}
By the mapping $u$ $=$ {\it{the constant}} on $\partial \mathbb{B}^{d}$,
 we thus deduce
  that the WHHF $u$ converges strongly to {\it{it}}
   as $t \nearrow \infty$ in $L^2(\mathbb{B}^d)$. \qed

%% file: ref_hhf_arch.tex
%
%#! platex hhf
%
%%%%%%%%%%%%%%%%%%%%%%%%%%%%%%%%%%%%%%%%%%%%%%%%%%%%%%%%%%%%%%%5
%
%\thanks{\it{I appreciate Prof.Xinfu Chen for discussion}}
\bibliographystyle{amsalpha}

\AD